\documentclass[11pt,leqno]{article}

\usepackage{graphicx}
\usepackage{latexsym} 
\usepackage{amsmath}
\usepackage{amssymb}
\usepackage{amsfonts}
\usepackage{enumerate}

\DeclareMathOperator{\sech}{sech}

\usepackage{theorem}

\theoremstyle{plain}
\newtheorem{teo}{Theorem}[section]
\newtheorem{lem}[teo]{Lemma}
\newtheorem{cor}[teo]{Corollary}
\newtheorem{prop}[teo]{Proposition}

{\theorembodyfont{\rmfamily}

}

\renewcommand{\eqref}[1]{\textnormal{(\ref{#1})}}

\numberwithin{equation}{section}

\newcommand{\cvd}{\hfill$\square$}
\newcommand{\proof}[1]{\noindent\textsc{Proof#1}}
\newcommand{\rmi}{\mathrm{i}}
\newcommand{\rme}{\mathrm{e}}
\newcommand{\rmd}{\mathrm{d}}

\title{Stable determination of a scattered wave from its far-field pattern: the high frequency asymptotics}

\author{Luca \textsc{Rondi}\thanks{Dipartimento di Matematica e Geoscienze,
Universit\`a degli Studi di Trieste, via Valerio, 12/1, 34127
Trieste, Italy. E-mail: \texttt{rondi@units.it}} \and
Mourad \textsc{Sini}\thanks{Johann Radon Institute for Computational and Applied Mathematics (RICAM),
 Austrian Academy of Sciences, Altenbergstr. 69, A-4040 
 Linz, Austria. E-mail: \texttt{mourad.sini@oeaw.ac.at}}}

\date{}

\begin{document}

\maketitle

\setcounter{section}{0}
\setcounter{secnumdepth}{2}

\begin{abstract}

We deal with the stability issue for the determination of outgoing time-harmonic acoustic waves from their far-field patterns. We are especially interested in keeping as explicit as possible the dependence of our stability estimates on the wavenumber of the corresponding Helmholtz equation and in understanding the high wavenumber, that is frequency, asymptotics.

Applications include stability results for the determination from far-field data of solutions of direct scattering problems with sound-soft obstacles and an instability analysis for the corresponding inverse obstacle problem. 

The key tool consists of establishing precise estimates on the behavior of Hankel functions with large argument or order.

\medskip

\noindent\textbf{AMS 2000 Mathematics Subject Classification}
Primary 35P25. Secondary 35R30.

\medskip

\noindent \textbf{Keywords} Helmholtz equation, outgoing solution, far-field, scattering problems,  stability estimates, high frequency, inverse scattering, Hankel functions.
\end{abstract}

\section{Introduction}\label{intro}

In recent years there has been an increasing attention to the study of how stability estimates for ill-posed problems involving the Helmholtz equation or the Schr\"odinger equation improve as the wavenumber and frequency or the energy, respectively, grows and might become extremely large.

One of the first rigorous justification of this phenomenon is due to Isakov and collaborators, \cite{Hry-Isa,Sub-Isa}, and concerns the Cauchy problem for the Helmholtz equation.

For what concerns corresponding inverse problems, increasing stability properties has been shown in many cases by many different authors. For instance, in
\cite{Bao-Lin-Tri}, an inverse source problem for the Helmholtz equation in the high frequencies regime was studied.
The inverse problem of determining the potential in a Schr\"odinger equation by boundary data 
in the high energies regime was considered in \cite{Isak11,Isak-Nag-Uhl-Wan,Isak-Wan} by geometrical optics techniques and, with a different method, in \cite{Isaev-Nov}. Let us notice that Isaev \cite{Isaev} developed a corresponding instability analysis showing the optimality of the previous estimates.
In \cite{Nag-Uhl-Wan} similar stability estimates were proved for the inverse problem of determining by boundary data an inhomogeneous medium for the acoustic wave equation.

We recall that typically these inverse problems are exponentially ill-posed and this is one of the main difficulties for numerical reconstruction. Such an instability character of these kinds of inverse problems was proved by Mandache, \cite{Man}, for the the Schr\"odinger equation at zero energy and the inverse conductivity problem, see also \cite{DC-R1} for other elliptic inverse boundary value problems and inverse scattering problems.
This is the main reason for trying to obtain a more stable reconstruction by changing the frequency or energy and in particular by using
high frequencies or energies. This motivated the search for  
stability estimates with an explicit dependence on the frequency or energy and for their 
high frequency or energy asymptotics. It has been shown that, asymptotically,
stability estimates may be expressed as the sum of a stable H\"older or Lipschitz term with a logarithmic one which is converging to zero as the frequency or energy tends to infinity, thus making the exponential ill-posedness less severe.

Another motivation can be found in the so-called multifrequency methods (called also hopping type algorithms) to reconstruct materials or interfaces from their scattered or far-field patterns. 
The main advantage of using such multifrequency data is that it can help to obtain accurate reconstructions without the need for a good initial guess. 
Different reconstruction methods using multifrequency data have been proposed in the last two decades or so, see for instance \cite{C-L:IEEE1995, Chen:IP1997,Bao-Hou-Li, B-T:JCM2010, Bao-Lin-Tri:2011,  S-T:IPI2012}. 
The convergence analysis of this type of algorithms was investigated 
in \cite{B-T:JCM2010, S-T:IPI2012} for the so-called recursive linearization algorithm proposed in \cite{Chen:IP1997}. In the analysis of these methods, 
the need for lower estimating the singular values of the linearized scattering problem in terms of the used frequencies arises naturally, see also \cite{S-T:2013}. 
In order to derive such a lower estimate, a crucial step is to estimate the scattered waves
from their far-field patterns, up to the boundary of the scatterers and with an explicit dependence upon the frequency.

In this paper we are mainly interested in the direct and inverse acoustic scattering problems for impenetrable scatterers, in particular sound-soft obstacles. For the inverse problem numerical evidence that the stability improves as the frequency grows was shown in
\cite{Col-Had-Pia}. It would be extremely interesting to rigorously prove such a phenomenon by establishing suitable stability estimates in the high frequencies regime. Unfortunately this seems to be still a challenging open problem. Nevertheless, for the direct scattering problem we obtain several interesting results. Our aim is to obtain stability estimates for the values of scattered waves from their far-field patterns, in the high frequency case and on the whole exterior of the scatterer. We are able to establish this result at least for smooth convex obstacles.

A crucial step, and one of the main results of the paper, is 
to prove stability estimates depending on the frequency for the determination of the near-field of an outgoing acoustic time-harmonic wave by its corresponding far-field. Such a problem has been solved for a fixed frequency by Isakov, \cite{Isak92}, see also \cite{Bus}. Very recently and independently Isakov \cite{Isak13} considered the high frequency case and showed that the stability improves as the frequency increases, even if his result is weaker than ours. In fact our estimates improve exponentially with respect to $k$, rather than polynomially. Besides, our a priori bound is of $L^2$ type instead of an $H^1$ type.
Moreover, in the regime where the Isakov's estimate is more meaningful, that is when $k$ is higher
than $\log(1/\varepsilon)$, $\varepsilon$ being an estimate of the norm of the far-field pattern,
we are able to obtain a Lipschitz stability estimate completely independent of $k$, see for instance Corollary~\ref{Lipschitz}. Let us finally observe that the stability estimate of Isakov does not show any improvement as $k$ grows if his a priori bound is of order higher than $\sqrt{k}$, which is often the case as we shall show in Section~\ref{uptotheboundarysec}. On the contrary, using the exponential improvement with respect to $k$, we are still able to obtain a Lipschitz stability result provided our a priori bound grows polynomially with respect to $k$, a fact that we shall
prove to hold at least for scattering solutions related to sound-soft obstacles, see
Section~\ref{uptotheboundarysec}.

Finally, by an instability analysis, we are able to evaluate from below how much the stability for the inverse scattering problem could improve as the frequency increases.

One of the main features of the paper is that
in all our results we keep the dependence on the frequency as explicit as possible. We also establish estimates for the full range of frequencies, with a particular attention to the case of high frequencies. We notice that most of the results present three different regimes. The regime of bounded frequencies where the usual ill-posedness shows up. A regime of high frequencies, with a limited improvement in the stability results, and a regime of extremely high frequencies where the improvement is much more significant.
Furthermore, most of the results are written for any space dimension $N\geq 2.$

Let us describe in more details the results of the paper.
We consider time-harmonic acoustic waves in a homogeneous and isotropic medium in a subset of $\mathbb{R}^N$, $N\geq 2$. Such a wave is characterized by its corresponding field $u$ which solves the \emph{reduced wave equation} or \emph{Helmholtz equation}
$$\Delta u+k^2u =0 $$
where $k>0$ is the \emph{wavenumber}. We recall that the wavenumber $k$ is the ratio between the corresponding frequency and the speed of sound.

We say that a time-harmonic acoustic wave in the exterior of a ball is \emph{outgoing} if its corresponding field $u$ satisfies the so-called \emph{Sommerfeld radiation condition}, that is
$$\lim_{r\to\infty}r^{(N-1)/2}\left(\displaystyle{\frac{\partial u^s}{\partial r}}-
\rmi ku^s\right)=0,\quad r=\|x\|$$
where the limit is intended to hold uniformly for all directions $\hat{x}=x/\|x\|\in \mathbb{S}^{N-1}$.

The Sommerfeld radiation
condition allows to characterize the
asymptotic behavior of the outgoing acoustic wave, namely we have that
$$
u(x)=\frac{\rme^{\rmi k\|x\|}}{\|x\|^{(N-1)/2}}\left\{
u_{\infty}(\hat{x})+O\left(\frac{1}{\|x\|}\right)
\right\},
$$
as $\|x\|$ goes to $+\infty$, uniformly in all directions $\hat{x}=x/\|x\|\in \mathbb{S}^{N-1}$.
The function $u_{\infty}$ is defined on $\mathbb{S}^{N-1}$ and is referred to as the
\emph{far-field pattern} of the field $u$.

Outgoing waves play a fundamental role in acoustic scattering theory. In fact, let us assume that in a homogeneous and isotropic medium in $\mathbb{R}^N$ there is 
a \emph{scatterer} $\Sigma$, that is a compact subset of $\mathbb{R}^N$ such that $G=\mathbb{R}^N\backslash\Sigma$ is connected. We recall that a scatterer $\Sigma$ is said to be an \emph{obstacle} if $\Sigma$ is the closure of an open set.

Let us assume that we send through the medium a time-harmonic acoustic wave, which is called incident wave. We call $k>0$ its \emph{wavenumber} and $u^i$ its corresponding field, the \emph{incident field}. Usually the incident wave is either a point source wave or a plane wave. We shall mainly focus on the latter case. If  $\omega\in\mathbb{S}^{N-1}$ is the \emph{direction of propagation} of the time-harmonic plane wave, then the incident field is given by 
$u^i(x)=\rme^{\rmi k\omega\cdot x}$, $x\in\mathbb{R}^N$.

The presence of the scatterer perturbs the incident wave by producing a so-called scattered wave which is characterized by being an outgoing time-harmonic acoustic wave. Its field $u^s$ is called the \emph{scattered field}. The total wave is the superposition of the incident wave and the scattered wave, that is its corresponding field, the \emph{total field} $u$, is simply the sum of the incident field and the scattered field. Namely the total field satisfies the following
\begin{equation}\label{dirpbm}
\left\{\begin{array}{ll}
\Delta u+k^2u=0&\text{in }\mathbb{R}^N\backslash \Sigma\\
u=u^i+u^s&\text{in }\mathbb{R}^N\backslash \Sigma\\
B.C. &\text{on }\partial \Sigma\\
\lim_{r\to\infty}r^{(N-1)/2}\left(\displaystyle{\frac{\partial u^s}{\partial r}}-
\rmi ku^s\right)=0&r=\|x\|.
\end{array}\right.
\end{equation}
The boundary condition on $\partial\Sigma$ depends on the nature of the scatterer.
In this paper we mainly focus on the case of impenetrable \emph{sound-soft} scatterers to which corresponds the following homogeneous Dirichlet boundary conditions
$$u=0\quad\text{on }\partial\Sigma.$$
However other boundary conditions may appear in the applications,
for instance the homogeneous Neumann condition for impenetrable sound-hard scatterers, the more general impedance boundary condition for impenetrable scatterers or
transmission conditions for penetrable scatterers.

Since the scattered wave is outgoing the asymptotic behavior of the scattered field $u^s$ is determined by its far-field pattern $u^s_{\infty}$.

Let us recall that \eqref{dirpbm} is referred to as the direct scattering problem. The corresponding
inverse scattering problem consists of the determination of the position and shape of a scatterer given the far-field patterns of the scattered waves corresponding to one or more incident planar waves. In this paper, precisely in Section~\ref{instabilitysec}, we treat the case of sound-soft scatterers and in order to perform more scattering, that is far-field, measurements,
we let vary the direction of propagation of the incident fields, keeping the same wavenumber $k>0$.
Let us notice that, instead of far-field measurements, one may also use so-called near-field measurements, that is the Cauchy data
of the scattered fields $u^s$ on the boundary of a domain containing the scatterer $\Sigma$ or, equivalently, the values of the scattered fields on a neighborhood of the boundary of such a domain.

One of the main results of the paper is a stability estimate for the determination of near-field data from far-field data for any outgoing solution of the Helmholtz equation, see Section~\ref{stabsec}. The technique used is the one developed by Isakov in \cite{Isak92}, see also \cite{Bus}. The main novelty here, besides the fact that our results hold for any space dimension $N\geq 2$, is that we investigate how the estimate changes with respect to the wavenumber $k$. In Theorem~\ref{Isakov} we deal with the case of $k$ belonging to a fixed compact interval of positive numbers.
This result is essentially a rephrasing of Isakov's result extended to any dimension $N\geq 2$ and in fact we obtain the usual logarithmic type estimate.
Then we deal with the high frequency, or wavenumber, case and we notice that the stability estimate improves as $k$ increases. Actually there are two regimes: the high frequencies regimes and the extremely high frequencies one. The first regime is treated in Theorem~\ref{mainteo} and it holds for wavenumbers $k$ which are at most of the order of $\log(1/\varepsilon)$, $\varepsilon$ being an estimate of the norm of the far-field pattern. Here the estimate is still of logarithmic type and improves as $k$ increases in an exponential way with respect to $k$. For the largest value of $k$ for which such a regime holds, the improvement leads to a H\"older type estimate. If $k$ is beyond such a threshold, that is in the extremely high frequencies regime, the stability estimate improves even further, see Proposition~\ref{highfreqprop}. In fact it is still at least of H\"older type and it may be actually written as the sum of a Lipschitz term plus one which is exponentially decaying with respect to $k$. Finally, if the a priori bound is of the order of some power of $k$, and $k$ is at least of the order of $\log(1/\varepsilon)$, we obtain a Lipschitz estimate fully independent of $k$, see Corollary~\ref{Lipschitz}.

Let us recall the basic idea of Isakov for obtaining this kind of stability estimates.
The key tool is using separation of variables for describing outgoing solutions of the Helmholtz equation. The dependence on the radial variable $r$ is given through suitable Hankel functions of first kind evaluated at $kr$. If the wavenumber is below a fixed constant, the stability estimates follows by studying the asymptotic behavior of these Hankel functions with respect to their order, a classical result in the theory of special functions. However, if we let the wavenumber tend to infinity the analysis is much more involved since we need to consider the asymptotic behavior of Hankel functions in three different regimes. In the first the argument, that is $kr$, is much larger than the order. In the second the argument and the order are both large but of the same magnitude. In the third the order is much larger than the argument. Such an asymptotic analysis is performed in Section~\ref{Hankel} whose main result is Theorem~\ref{regimes} which contains the asymptotic behavior of Hankel functions in the previous three regimes. This is the main technical result of the paper and the essential tool for obtaining the estimates of Section~\ref{stabsec}.

In Section~\ref{uptotheboundarysec} we apply the previous estimates to direct scattering problems. We assume that $\Sigma_1$ and $\Sigma_2$ are two sound-soft obstacles and that $K$ is the convex hull of their union. We assume that $\Sigma_1$ and $\Sigma_2$ are star-shaped and smooth enough.
Let $u_1$ and $u_2$ be the solutions to the direct scattering problem with $\Sigma$ replaced by $\Sigma_1$ and $\Sigma_2$, respectively. Our aim is to estimate the difference of $u_1$ and $u_2$ up to the boundary of $K$ by the difference of the
far-field patterns of the corresponding scattered waves.
We use the results of Section~\ref{stabsec} to proceed from far-field to near-field data and then we use the results of \cite{Hry-Isa,Sub-Isa} to estimate, from the the near-field data, $u_1-u_2$ up to the boundary of $K$. We obtain a stability estimate with an explicit dependence on the wavenumber $k$. However we need to note that, for the time being, we are not able to prove any increasing stability property as $k$ grows and tend to infinity, see Theorem~\ref{uptheboundaryteo} and the following discussion.

About our stability estimate, we wish to use as less as possible a priori information, namely only a priori bounds of $L^2$ type on the solutions. For this reason we use an integral norm to estimate the difference between $u_1$ and $u_2$. In order to use the results of \cite{Hry-Isa,Sub-Isa} on the whole exterior of $K$ some technical difficulties arise which are solved by using the $L^1$ norm, instead of the $L^2$ norm, to estimate the difference between $u_1$ and $u_2$ and by the help of a technical geometrical lemma, Lemma~\ref{Kproperties}.

It would be desirable to proceed further with the analysis and obtain suitable stability estimates, with a precise dependence on $k$, up to the boundary of the unbounded connected component of
$\mathbb{R}^N\backslash(\Sigma_1\cup\Sigma_2)$ or
for the corresponding inverse scattering problem.
However, both seem to be rather difficult open problems.

Let us notice that a rather long preliminary part of this section contains the essential a priori estimates which are needed to implement the previously described strategy. We believe that this part may also be of independent interest. It is here that the main assumptions on $\Sigma_1$ and $\Sigma_2$, namely star-shapedness and smoothness, are needed. The key ingredient in the high frequencies regime is provided by Theorem~\ref{Cha-Monteo} which follows from results due to Chandler-Wilde and Monk, \cite{Cha-Mon}, and to Melenk, \cite{Mel}.

Finally, on Section~\ref{instabilitysec} we perform an instability analysis for the corresponding inverse scattering problem with sound-soft obstacles, extending the instability result in \cite{DC-R1} to the high frequencies regime. Let us recall that a similar result has been obtained by Isaev for the inverse problem of  determining by boundary data the potential in 
the Schr\"odinger equation in the high energies regime, see \cite{Isaev}.

The main results of this section, Theorem~\ref{instmainteo} and Corollary~\ref{instcorollary}, show that the instability improves as the wavenumber $k$ increases. The high frequencies regime holds for wavenumbers $k$ which are at most of the order of $\log(1/\varepsilon)$, $\varepsilon$ being an estimate of the norm of the error in the far-field pattern for all possible directions of propagation of the incident field. In this regime the improvement is not so significant, however, beyond it, that is for extremely high frequencies, the improvement is more relevant since the logarithmic instability term has a multiplicative constant converging to zero, as $k$ goes to $+\infty$, in a polynomial way with respect to $k$.

Let us notice that the norm used to estimate the error in the far-field pattern is an arbitrary $H^s$ Sobolev norm, with $s\geq 0$. For $s$ sufficiently large this is stronger than the $L^{\infty}$ norm, that is our result applies also when $\varepsilon$ is an estimate of the superior, for all 
directions of propagation of the incident field, of
the error in the far-field pattern measured in the $L^{\infty}$ norm. Moreover, we remark that our instability results hold for star-shaped and even convex obstacles. Even if the reconstruction of star-shaped or convex obstacle is considered to be more stable, our analysis shows that star-shapedness or convexity do not provide a significant advantage.

The proof is based on the original idea by Mandache, \cite{Man}, which has been generalized and applied to inverse scattering problems in \cite{DC-R1}. However again we need to take into account the fact that the wavenumber may be arbitrarily large, therefore a careful use of the results of Section~\ref{Hankel} is required.

The plan of the paper is the following. We begin with a preliminary section, Section~\ref{farfield},
where we give a separation of variables description of scattered waves and their far-field patterns. Moreover we introduce suitable classes of smooth star-shaped obstacles. In Section~\ref{Hankel}
we develop the asymptotic analysis for Hankel functions. The stability estimates for the determination of near-field data from far-field data are contained in Section~\ref{stabsec}. In Section~\ref{uptotheboundarysec} we consider the application to direct scattering problems with sound-soft scatterers. In particular we apply the estimate of the previous section to the stable determination of scattered waves from their far-field patterns. Finally, in Section~\ref{instabilitysec}
the instability analysis of the corresponding inverse scattering problem is considered.

\subsubsection*{Acknowledgements}
L.R. is partly supported by Universit\`a degli Studi di Trieste through Fondo per la Ricerca di Ateneo -- FRA 2012.
M.S. is partly supported by the Johann Radon Institute for Computational and Applied Mathematics (RICAM) of the 
 Austrian Academy of Sciences.
Part of this work
was done while L.R. was visiting RICAM, whose support and hospitality is gratefully acknowledged.

\section{Preliminaries}\label{farfield}

Let us fix an integer $N\geq 2$. For any $x\in\mathbb{R}^N$ and any $s>0$, 
$B_s(x)$ denotes the ball contained in $\mathbb{R}^N$ with radius $s$ and center $x$.
Moreover, $B_s=B_s(0)$ and, finally, for any $E\subset \mathbb{R}^N$, we denote $B_s(E)=\bigcup_{x\in E}B_s(x)$. For any $E\subset \mathbb{R}^N$, $|E|$ denotes as usual the $N$-dimensional Lebesgue measure of $E$.

We fix a \emph{scatterer} $\Sigma$ in $\mathbb{R}^N$, that is a compact subset of $\mathbb{R}^N$ such that $G=\mathbb{R}^N\backslash\Sigma$ is connected. We recall that a scatterer $\Sigma$ is said to be an \emph{obstacle} if $\Sigma$ is the closure of an open set. We fix the \emph{wavenumber} $k>0$ and a \emph{direction of propagation} $\omega\in\mathbb{S}^{N-1}$. Then the \emph{incident field} is
the time-harmonic plane wave $u^i(x)=\rme^{\rmi k\omega\cdot x}$, $x\in\mathbb{R}^N$.
The \emph{total field} $u=u(\omega,k,\Sigma)$ is the sum of the incident field and of the \emph{scattered field} $u^s$ and is the solution to the following exterior boundary value problem
\begin{equation}\label{Helmeq}
\left\{\begin{array}{ll}
\Delta u+k^2u=0&\text{in }\mathbb{R}^N\backslash \Sigma\\
u=u^i+u^s&\text{in }\mathbb{R}^N\backslash \Sigma\\
u=0&\text{on }\partial \Sigma\\
\lim_{r\to\infty}r^{(N-1)/2}\left(\displaystyle{\frac{\partial u^s}{\partial r}}-
\rmi ku^s\right)=0&r=\|x\|
\end{array}\right.
\end{equation}
where the last limit, the so-called \emph{Sommerfeld radiation condition},
holds
uniformly for all directions $\hat{x}=x/\|x\|\in \mathbb{S}^{N-1}$.
We remark that the homogeneous Dirichlet boundary condition corresponds to a \emph{sound-soft}
scatterer $\Sigma$.

The Sommerfeld radiation
condition characterizes outgoing waves and implies that the
asymptotic behavior of the scattered field is given by
\begin{equation}\label{asympt}
u^s(x;\omega,k,\Sigma)=\frac{\rme^{\rmi k\|x\|}}{\|x\|^{(N-1)/2}}\left\{
u^s_{\infty}(\hat{x};\omega,k,\Sigma)+O\left(\frac{1}{\|x\|}\right)
\right\},
\end{equation}
as $\|x\|$ goes to $+\infty$, uniformly in all directions $\hat{x}=x/\|x\|\in \mathbb{S}^{N-1}$.
The function $u^s_{\infty}$ is called the \emph{far-field pattern} of the scattered field $u^s$ of
the solution to \eqref{Helmeq}.

For any sound-soft scatterer $\Sigma$,
we denote
$\mathcal{A}(\Sigma):\mathbb{S}^{N-1}\times \mathbb{S}^{N-1}\times (0,\infty)\mapsto\mathbb{C}$ its
\emph{far-field map}, that is, for any $\hat{x}$,
$\omega\in \mathbb{S}^{N-1}$ and any $k>0$,
\begin{equation}\label{softffpdef}
\mathcal{A}(\Sigma)(\hat{x},\omega,k)=u^s_{\infty}(\hat{x};\omega,k,\Sigma),
\end{equation}
where $u^s_{\infty}$ is the far-field pattern of the scattered field $u^s$ of
the solution to \eqref{Helmeq}.

Let us remark that the following \emph{reciprocity relation} holds, see for instance
\cite[Theorem~3.13]{Col e Kre98}. For any scatterer $\Sigma$ and any $k\in(0,\infty)$ we
have
\begin{equation}\label{recrel}
\mathcal{A}(\Sigma)(\hat{x},\omega,k)=\mathcal{A}(\Sigma)(-\omega,-\hat{x},k)\quad\text{for
any }\hat{x},\omega\in \mathbb{S}^{N-1}.
\end{equation}

Moreover, we have the following characterization of the $L^2$ norm of the far-field pattern, see for instance \cite[Theorem~3.2.1]{Ned} for $N=3$,
\begin{equation}\label{Nedelec}
\|\mathcal{A}(\Sigma)(\cdot,\omega,k)\|^2_{L^2(\mathbb{S}^{N-1})}=2\left(\frac{2\pi}{k}\right)^{(N-1)/2}
\Im\left(
\rme^{(N-3)\pi\rmi/4}
\mathcal{A}(\Sigma)(\omega,\omega,k)
\right).
\end{equation}

We wish to decompose the far-field pattern in spherical harmonics. Fixed $N\geq 2$, let us consider the orthonormal basis of $L^2(\mathbb{S}^{N-1})$
\begin{equation}\label{spherarm}
\{v_i\}_{i\in\mathbb{N}}=\{f_{jp}:\ j\geq 0\text{ and }1\leq p \leq p_j\}
\end{equation}
that consists of (real-valued) spherical harmonics, that is
each $f_{jp}$ is a \emph{spherical harmonic} of degree $j$, $j$ being a
nonnegative integer, such that
$\|f_{jp}\|_{L^2(\mathbb{S}^{N-1})}=1$. The elements $v_i$, $i\in\mathbb{N}$, are ordered in the natural way.

The integers $p_j$
are the dimensions of the spaces of spherical harmonics of degree $j$ and we have that,
see for instance \cite[page~4]{Mul},
$$p_j=\left\{\begin{array}{ll}1 & \text{if }j=0,\\
\frac{(2j+N-2)(j+N-3)!}{j!(N-2)!} & \text{if }j\geq 1,\end{array}\right.$$
so that
$$p_j\leq 2(j+1)^{N-2},\quad j\geq 0,$$
and
\begin{equation}\label{numspherarm}
\sum_{j=0}^np_j\leq\sum_{j=0}^n2(j+1)^{N-2}\leq 2(n+1)^{N-1},\quad\text{for
any }n\in\mathbb{N}.
\end{equation}

For any spherical harmonic $f$, we call $\gamma(f)$ the degree of the spherical harmonic, that is
$\gamma(f_{jp})=j$.
We have that $\gamma(v_i)$ is an increasing sequence, with
respect to $i$, whose asymptotic behavior satisfies the following property.
Fixed $n\in\mathbb{N}$, we have that, by \eqref{numspherarm},
$\#\{i\in\mathbb{N}:\ \gamma(v_i)\leq n\}$ is clearly bounded from
above by $2(n+1)^{N-1}$, therefore
\begin{equation}\label{polygrowth}
\#\{i\in\mathbb{N}:\ \gamma(v_i)\leq n\}\leq C(n+1)^{p}
\end{equation}
with $C=2$ and $p=N-1$. We recall that $\#$ denotes the number of elements.

We recall that the function
\begin{equation}\label{polyharm}
u_{jp}(x)=\|x\|^jf_{jp}(x/\|x\|)
\end{equation}
is harmonic in $\mathbb{R}^N$ and solves the following eigenvalue problem in $B_1$
\begin{equation}
\Delta u_{jp}=0\text{ in }B_1;\quad
\frac{\partial u_{jp}}{\partial \nu}=ju_{jp}\text{ on }\partial B_1.
\end{equation}

For any function $g$ belonging to $L^2(\mathbb{S}^{N-1}\times \mathbb{S}^{N-1})$,
we decompose it in spherical harmonics in the following way
$$
g(\hat{x},\omega)=\sum_{i,l}a_{i,l}v_i(\hat{x})v_l(\omega)\quad (\hat{x},\omega)\in  \mathbb{S}^{N-1}\times \mathbb{S}^{N-1}
$$
where the complex-valued coefficients $a_{i,l}$ are given
by
$$
a_{i,l}=\int\!\!\!\int_{\mathbb{S}^{N-1}\times \mathbb{S}^{N-1}}g(\hat{x},\omega)
v_i(\hat{x})v_l(\omega)\mathrm{d}\hat{x}
\mathrm{d}\omega.
$$
For any $s\geq 0$ we define the norm of the Sobolev space $H^s(\mathbb{S}^{N-1}\times \mathbb{S}^{N-1})$ as follows
\begin{equation}\label{Hsnorm}
\|g\|^2_{H^s(\mathbb{S}^{N-1}\times \mathbb{S}^{N-1})}=\sum_{i,l}(1+\gamma(v_i)+\gamma(v_l))^{2s}|a_{i,l}|^2.
\end{equation}
Then for any $s\geq 0$ we call $Y_s(\mathbb{S}^{N-1}\times \mathbb{S}^{N-1})$ the space
$$Y_s(\mathbb{S}^{N-1}\times \mathbb{S}^{N-1})=\{g\in L^2(\mathbb{S}^{N-1}\times \mathbb{S}^{N-1}):\ \|g\|_{s}<+\infty\}$$ 
where
$$\|g\|_{s}=\sup_{i,l} \left((1+\max\{\gamma(v_i),\gamma(v_l)\})^{2s+N-1/2}|a_{i,l}|\right)$$
We notice that $Y_s(\mathbb{S}^{N-1}\times \mathbb{S}^{N-1})\subset H^s(\mathbb{S}^{N-1}\times \mathbb{S}^{N-1})$
and the immersion is continuous, in fact we have
$$\|g\|_{H^s(\mathbb{S}^{N-1}\times \mathbb{S}^{N-1})}\leq 
4\|g\|_{s}\quad
\text{for any }g\in Y_s(\mathbb{S}^{N-1}\times \mathbb{S}^{N-1}),
$$
since
\begin{multline*}
\sum_{i,l}(1+\gamma(v_i)+\gamma(v_l))^{2s}|a_{i,l}|^2\leq
\sum_{i,l}(1+\max\{\gamma(v_i),\gamma(v_l)\})^{4s}|a_{i,l}|^2\leq\\
16\sup_{i,l}\left((1+\max\{\gamma(v_i),\gamma(v_l)\})^{4s+2N-1}|a_{i,l}|^2\right)
\end{multline*}
see \cite[page~1439]{Man}.

The decomposition of the far-field pattern in spherical harmonics is given by,
for any $(\hat{x},\omega,k)\in \mathbb{S}^{N-1}\times \mathbb{S}^{N-1}\times (0,\infty)$,
\begin{equation}\label{decomposition}
\mathcal{A}(\Sigma)(\hat{x},\omega,k)=\sum_{i,l}b_{i,l}(k)v_i(\hat{x})v_l(\omega),
\end{equation}
where the complex-valued coefficients $b_{i,l}(k)$ are given, for any $k\in (0,\infty)$,
by
\begin{equation}\label{decompcoeff}
b_{i,l}(k)=\int\!\!\!\int_{\mathbb{S}^{N-1}\times \mathbb{S}^{N-1}}\mathcal{A}(\Sigma)(\hat{x},\omega,k)
v_i(\hat{x})v_l(\omega)\mathrm{d}\hat{x}
\mathrm{d}\omega.
\end{equation}
Furthermore, we use the following characterization
\begin{equation}\label{decompcoeff2}
b_{i,l}(k)=\int_{\mathbb{S}^{N-1}}\tilde{b}_i(\omega,k)
v_l(\omega)\mathrm{d}\omega,
\end{equation}
where
the complex-valued coefficients $\tilde{b}_i(\omega,k)$ are, for any
$\omega\in \mathbb{S}^{N-1}$ and any
$k\in (0,\infty)$, the Fourier coefficients, with respect to the orthonormal basis
$\{v_i\}_{i\in\mathbb{N}}$,
of the far-field pattern $u^s_{\infty}(\cdot;\omega,k,\Sigma)$ corresponding to
the scattered field of the solution to \eqref{Helmeq}, that is
\begin{equation}\label{decompcoeff3}
\tilde{b}_i(\omega,k)=\int_{\mathbb{S}^{N-1}}\mathcal{A}(\Sigma)(\hat{x},\omega,k)
v_i(\hat{x})\mathrm{d}\hat{x}.
\end{equation}

Let us assume that $\Sigma\subset\overline{B_R}$ for some positive constant $R$.
Then for any $x\in\mathbb{R}^N\backslash \overline{B_R}$ we have
\begin{equation}\label{decompscat}
u^s(x;\omega,k,\Sigma)=
\sum_i\hat{b}_i(\omega,k)\frac{H^{(1)}_{\gamma(v_i)+(N-2)/2}(k\|x\|)}{(k\|x\|)^{(N-2)/2}}
v_i(x/\|x\|),
\end{equation}
where $\hat{b}_i=\hat{b}_i(\omega,k)$ are complex-valued coefficients given by
\begin{equation}\label{decompscatcoeff}
\hat{b}_i(\omega,k)
\frac{H^{(1)}_{\gamma(v_i)+(N-2)/2}(kr)}{(kr)^{(N-2)/2}}
=\int_{\mathbb{S}^{N-1}}u^s(r\hat{x};\omega,k,\Sigma)
v_i(\hat{x})\mathrm{d}\hat{x},\quad\text{for any }r>R,
\end{equation}
where, for any real $\nu\geq 0$, $H^{(1)}_\nu$ denotes the \emph{Hankel function} of
first kind and order $\nu$.

The relationship between
coefficients $\tilde{b}_i$ and $\hat{b}_i$ is the following
\begin{equation}\label{link}
\tilde{b}_i(\omega,k)=(\pi/2)^{-1/2}k^{-(N-1)/2}(-\rmi)^{\gamma(v_i)+(N-1)/2}
\hat{b}_i(\omega,k).
\end{equation}
Therefore, for any $r>R$,
\begin{equation}\label{decompscat2}
u^s(r\hat{x};\omega,k,\Sigma)=
\sum_i\tilde{b}_i(\omega,k)
(\pi k/2)^{1/2}\rmi^{\gamma(v_i)+(N-1)/2}
\frac{H^{(1)}_{\gamma(v_i)+(N-2)/2}(kr)}{r^{(N-2)/2}}
v_i(\hat{x}),
\end{equation}
hence
\begin{equation}\label{decompscat3}
\|u^s(\cdot;\omega,k,\Sigma)\|^2_{L^2(\partial B_r)}= \frac{\pi}{2}
\sum_i|\tilde{b}_i(\omega,k)|^2
kr
\left|H^{(1)}_{\gamma(v_i)+(N-2)/2}(kr)\right|^2.
\end{equation}
Obviously, by \eqref{softffpdef} and \eqref{decompcoeff3} we have
\begin{equation}\label{farfielderror}
\|u^s_{\infty}(\cdot;\omega,k,\Sigma)\|^2_{L^2(\mathbb{S}^{N-1})}=\sum_i|\tilde{b}_i(\omega,k)|^2.
\end{equation}

In order to estimate the near-field from the far-field, a crucial step is to estimate
the asymptotic behavior of Hankel functions $H^{(1)}_{\nu}(z)$ as the order  $\nu$ goes to infinity, for a fixed $z=kr$. On the other hand, we are interested in the corresponding estimate as the wavenumber $k$ goes to infinity, therefore we need to consider the behavior of 
$H^{(1)}_{\nu}(z)$ with the order $\nu$ and the argument $z=kr$ which may be both large at the same time. We deal with this issue in the following section.

We conclude this section by introducing suitable classes of obstacles.
We fix integers $N\geq 2$ and
$m\geq 1$ and positive constants $\beta$, $R_0$ and $\delta$, $R_0<R_0+\delta\leq \beta$.

Let $g$ be a strictly positive continuous function defined on
$\mathbb{S}^{N-1}$. Let $\Sigma(g)$ be the compact set given by the \emph{radial subgraph} of $g$, that is
$$\Sigma(g)=\{y\in\mathbb{R}^N:\
y=\rho\omega,\ 0\leq \rho\leq g(\omega),\ \omega\in
\mathbb{S}^{N-1}\}.$$

We denote
$$X(m,\beta,R_0,\delta)=
\{\Sigma(g):\ g\in C^m(\mathbb{S}^{N-1}),
\ \|g\|_{C^m(\mathbb{S}^{N-1})}\leq \beta
\text{ and }R_0\leq g\leq R_0+\delta\}.$$

We notice that $X=X(m,\beta,R_0,\delta)$ is a metric space, endowed with the Hausdorff distance $d_H$,
that consists of obstacles in $\mathbb{R}^N$ which are star-shaped with
respect to the origin. Moreover, for any integer $m\geq 3$, and for any positive $\beta$ and $R_0$,
there exists $\tilde{\delta}>0$, depending on $N$, $m$, $\beta$ and $R_0$ only, such that
if $0<\delta\leq\tilde{\delta}$, then 
any element of 
$X$ is even convex.

Let us also notice that there exists a constant $E$, depending on $N$, $m$, $\beta$ and $R_0$ only, such that
\begin{equation}\label{areabound}
\mathcal{H}^{N-1}(\Sigma)\leq E
\quad\text{for any }\Sigma\in X(m,\beta,R_0,\delta),
\end{equation}
where $\mathcal{H}^{N-1}$ denotes the $(N-1)$-dimensional Hausdorff measure.

Given the metric space $(X,d_H)$ and a positive $\delta$, $X'\subset X$ is said to be
$\delta$-\emph{discrete} if any two distinct points $x_1$, $x_2\in X'$ satisfy
$d_H(x_1,x_2)\geq\delta$. We have the following result, whose proof may be easily obtained following the arguments of the proof of Lemma~2 in \cite{Man}.

\begin{prop}\label{discretesetprop}
Let us fix integers $N\geq 2$ and
$m\geq 1$ and positive constants $\beta$ and $R_0$.

Then, there exists a positive constant $\delta_0$, depending on
$N$, $m$, $\beta$ and $R_0$
only, such that 
$R_0+\delta_0\leq \beta$ and
for any
$\delta$, $0<\delta\leq\delta_0$, we can find a $\delta$-discrete subset of
$X(m,\beta,R_0,\delta)$ with at least $\exp(2^{-N}\delta_0^{(N-1)/m}\delta^{-(N-1)/m})$ elements.
\end{prop}

\section{Estimates on the asymptotic behavior of Hankel functions}\label{Hankel}

For any $\nu\geq 0$ and any $z>0$, let $H^{(1)}_{\nu}(z)$ be the Hankel function of first kind of order $\nu$ and argument $z$. Concerning basic properties of Hankel and Bessel functions we refer to \cite{Leb} and especially \cite{Wat}.

It is well-known that the following asymptotic behaviors of $H^{(1)}_{\nu}(z)$ holds true.
First, fixed $\nu\geq 0$, we have that
\begin{equation}\label{argument}
H^{(1)}_{\nu}(z)=\left(\frac{2}{\pi z}\right)^{1/2}\rme^{\rmi(z-\nu\pi/2-\pi/4)}\left[1+O(z^{-1})\right]
\quad\text{as }z\to+\infty.
\end{equation}
On the other hand, fixed $z>0$, we have that
\begin{equation}\label{order}
H^{(1)}_{\nu}(z)\sim-\rmi\left(\frac{2}{\pi\nu}\right)^{1/2}\left(\frac{\rme z}{2\nu}\right)^{-\nu}
\quad\text{as }\nu\to+\infty.
\end{equation}
where $\sim$ means that the quotient between the two functions tends to $1$, as $\nu\to\infty$.

We need similar estimates on the asymptotic behavior of Hankel functions which hold uniformly for suitable intervals of $z$ and $\nu$ respectively. We obtain a full hierarchy of asymptotic behaviors depending on the relationship between $z$ and $\nu$ as one or both of these parameters go to $+\infty$.

We begin by recalling that
$$J_{\nu}(z)=\sum_{k=0}^{+\infty}\frac{(-1)^k(z/2)^{\nu+2k}}{\Gamma(k+1)\Gamma(k+\nu+1)}\quad z>0,\ \nu\in\mathbb{R}$$
where $J_{\nu}$ is the Bessel function of first kind of order $\nu$. Here $\Gamma$ denotes the Gamma function.
We also recall that $Y_{\nu}$, the Bessel function of second kind of order $\nu$, is given by
$$Y_{\nu}(z)=\frac{J_{\nu}(z)\cos(\nu\pi)-J_{-\nu}(z)}{\sin(\nu\pi)}\quad z>0,\ \nu\in\mathbb{R}\backslash\mathbb{Z},$$
whereas
$$Y_{n}(z)=\lim_{\nu\to n}Y_{\nu}(z)\quad z>0,\ n\in\mathbb{Z}.$$
Finally,
$$H^{(1)}_{\nu}(z)=J_{\nu}(z)+\rmi Y_{\nu}(z)\quad z>0,\ \nu\in\mathbb{R}.$$

Let us recall here the basic properties of the Gamma function $\Gamma$ that will be repeatedly used in the sequel. We note that $\Gamma(z)$ is defined for any $z\in \mathbb{C}$ such that $z\neq 0,-1,-2,\ldots$, and we have that $\Gamma(1)=1$. The most important property is the following
\begin{equation}\label{mainpropGamma}
\Gamma(z+1)=z\Gamma(z)\quad\text{for any }z\neq 0,-1,-2,\ldots,
\end{equation}
from which we deduce that,  for any $n=0,1,2,\ldots$, $n!=\Gamma(n+1)$. We also use that
\begin{equation}\label{secondmainpropGamma}
\Gamma(z)\Gamma(1-z)=\frac{\pi}{\sin(\pi z)}\quad\text{for any }z\not\in \mathbb{Z}.
\end{equation}
Finally the following version of the classical Stirling formula will be used
\begin{equation}\label{Stirling}
\Gamma(x)=\sqrt{2\pi}x^{x-1/2}\rme^{-x}[1+r(x)],\quad x>0
\end{equation}
where
$$
|r(x)|\leq \rme^{1/(12 x)}-1.
$$

We begin with the following easy remarks. For any $a$, $0\leq a<1$, any $z>0$ and any $\nu>0$,
let
$$A=A(a,z,\nu)=\sum_{k=0}^{+\infty}\frac{(-1)^k(z^2/4)^k}{\Gamma(k+a+1)\Gamma(k+\nu+1)}.$$

Let us fix a constant $z_2$, $z_2\geq 1$. For any $\varepsilon$, $0<\varepsilon\leq 1/\rme$, using the basic properties of the $\Gamma$ function, we can find $\nu_0=\nu_0(\varepsilon,z_2)>0$ such that
for any $0<z\leq z_2$ and any $\nu\geq\nu_0$ we have
\begin{equation}\label{Aestimate}
(1-\varepsilon)\frac{1}{\Gamma(\nu+1)\Gamma(a+1)}\leq A(a,z,\nu)\leq
\frac{1}{\Gamma(\nu+1)\Gamma(a+1)}.
\end{equation} 
Notice that it is enough that $z_2^2/4\leq \varepsilon (\nu_0+1)$, that is for instance we may set
$$\nu_0(\varepsilon,z_2)=\frac{z_2^2}{2\varepsilon}.$$
We may immediately infer the following.

\begin{prop}\label{firstprop}
Let us fix $z_2\geq 1$ and $\varepsilon$, $0<\varepsilon\leq 1/\rme$.
Then there exists a positive constant $\tilde{\nu}_0=\tilde{\nu}_0(\varepsilon)$ such that, for any
$a$, $0\leq a<1$, 
for any $0<z\leq z_2$, and
for any $\nu\geq \max\{\nu_0(\varepsilon/2,z_2),\tilde{\nu}_0(\varepsilon)\}$, we have
\begin{equation}\label{Aestimate2}
\frac{(1-\varepsilon)}{\Gamma(a+1)}\left(\frac{1}{\sqrt{2\pi\nu}}\nu^{-\nu}\rme^{\nu}\right)
\leq A(a,z,\nu)\leq
\frac{(1+\varepsilon)}{\Gamma(a+1)}\left(\frac{1}{\sqrt{2\pi\nu}}\nu^{-\nu}\rme^{\nu}\right)
\end{equation}
and, consequently, picking $a=0$,
\begin{equation}\label{Jestimate}
\frac{(1-\varepsilon)}{\sqrt{2\pi\nu}}\left(\frac{\rme z}{2\nu}\right)^{\nu}
\leq J_{\nu}(z)\leq\frac{(1+\varepsilon)}{\sqrt{2\pi\nu}}\left(\frac{\rme z}{2\nu}\right)^{\nu}.
\end{equation}
Notice that $\nu\geq \rme z$, therefore
\begin{equation}\label{Jestimate2}
0\leq J_{\nu}(z)\leq\frac{(1+\varepsilon)}{\sqrt{2\pi\nu}},
\end{equation}
hence there exists $\tilde{\nu}_1=\tilde{\nu}_1(\varepsilon)$, $\tilde{\nu}_1\geq \tilde{\nu}_0$, such that if $\nu$ is also greater than or equal to $\tilde{\nu}_1$ we have
\begin{equation}\label{Jestimate3}
0\leq J_{\nu}(z)\leq\varepsilon.
\end{equation}
\end{prop}

Let us now fix $\nu=n-1/2$ for some positive integer $n$. We obtain that
$\cos(\nu\pi)=0$ whereas $\sin(\nu\pi)=(-1)^{n+1}$. Therefore, for any $z>0$,
$Y_{\nu}(z)=(-1)^nJ_{-\nu}(z)$. We notice that
$$Y_{\nu}(z)=(-1)^n(z/2)^{-\nu}
\left[\sum_{k=0}^{\infty}\frac{(-1)^k(z^2/4)^k}{\Gamma(k+1)\Gamma(k-\nu+1)}
\right]=(-1)^n(z/2)^{-\nu}\left[A+B\right]
\quad z>0,$$
where
$$B=\sum_{k=0}^{n-1}\frac{(-1)^k(z^2/4)^k}{\Gamma(k+1)\Gamma(k-\nu+1)}\quad\text{and}\quad A=\sum_{k=n}^{+\infty}\frac{(-1)^k(z^2/4)^k}{\Gamma(k+1)\Gamma(k-\nu+1)}.$$
Let us notice that by a simple change of the index, we have that
$$A=(-1)^n(z^2/4)^nA(1/2,z,n),$$
therefore
$$Y_{\nu}(z)=(z/2)^{2n-\nu}A(1/2,z,n)+(-1)^n(z/2)^{-\nu}B=\tilde{A}+\tilde{B}.$$

By Proposition~\ref{firstprop}, for any $z$, $0<z\leq z_2$, and any $\nu \geq \max\{\nu_0(\varepsilon/2,z_2),\tilde{\nu}_0(\varepsilon)\}$, we have
\begin{equation}\label{Atildeestimate}
\frac{(1-\varepsilon)}{\Gamma(3/2)}\sqrt{\frac{z}{4\pi n}}
\left(\frac{\rme z}{2n}\right)^{n}
\leq \tilde{A}\leq
\frac{(1+\varepsilon)}{\Gamma(3/2)}\sqrt{\frac{z}{4\pi n}}
\left(\frac{\rme z}{2n}\right)^{n}.
\end{equation}
We infer that
\begin{equation}\label{Atildeestimate2}
0\leq \tilde{A}\leq
\frac{(1+\varepsilon)}{\Gamma(3/2)}\sqrt{\frac{1}{2\rme\pi}}
\left(\frac{\rme z}{2n}\right)^{n+1/2}.
\end{equation}
hence there exists $\tilde{\nu}_2=\tilde{\nu}_2(\varepsilon)$, $\tilde{\nu}_2\geq \tilde{\nu}_1$, such that if $\nu$ is also greater than or equal to $\tilde{\nu}_2$ we have
\begin{equation}\label{Atildeestimate3}
0\leq \tilde{A}\leq\varepsilon.
\end{equation}

We continue by evaluating $B$ and, correspondingly, $\tilde{B}$.
A simple computation, where we use in particular \eqref{mainpropGamma} and also \eqref{secondmainpropGamma}, leads to
\begin{multline*}
\tilde{B}=-(z/2)^{-\nu}\sum_{k=0}^{n-1}\left(\frac{(z^2/4)^k}{\Gamma(k+1)\Gamma(1/2)}\prod_{i=1}^{n-1-k}(i-1/2)\right)=
\\-(z/2)^{-\nu}\sum_{k=0}^{n-1}\left(\frac{(z^2/4)^k\Gamma(n-k-1/2)}{\Gamma(k+1)(\Gamma(1/2))^2}\right)=-\frac{(z/2)^{-\nu}}{\pi}\sum_{k=0}^{n-1}\left((z^2/4)^k\frac{\Gamma(n-k-1/2)}{\Gamma(k+1)}\right).
\end{multline*}
Taking $h=n-1-k$, we have that
$$\tilde{B}=-\frac{(z/2)^{-\nu}}{\pi}\Gamma(\nu)\left[1+\sum_{h=0}^{n-2}\left((z^2/4)^{n-1-h}\frac{\Gamma(h+1/2)}{\Gamma(\nu)\Gamma(n-h)}\right)
\right]=-\frac{(z/2)^{-\nu}}{\pi}\Gamma(\nu)\left[1+R\right].$$
A lengthy but straightforward computation would allow us to estimate the remainder $R$, by studying separately the sum with $h$ below $(n-2)/2$ and the one with $h$ above $(n-2)/2$. We obtain that
$R$ goes to zero as $\nu$ goes to $+\infty$. However, for our purposes, at this time it is enough to note that $R\geq 0$, therefore
$$\tilde{B}\leq -\frac{(z/2)^{-\nu}}{\pi}\Gamma(\nu)\leq 
-(1-\varepsilon)\left(\frac{2}{\pi\nu}\right)^{1/2}\left(\frac{\rme z}{2\nu}\right)^{-\nu}$$
provided $\nu\geq \tilde{\nu}_3(\varepsilon)\geq \tilde{\nu}_2(\varepsilon)$ as well. Up to changing $\tilde{\nu}_3=\tilde{\nu}_3(\varepsilon)$, we obtain that
$$Y_{\nu}(z)\leq
-(1-\varepsilon)\left(\frac{2}{\pi\nu}\right)^{1/2}\left(\frac{\rme z}{2\nu}\right)^{-\nu}
$$
and, using \eqref{Jestimate3}, also that
$$|J_{\nu}(z)|\leq \varepsilon |Y_{\nu}(z)|.$$

Let us recall the following lemma, see for instance \cite[Appendix~B]{Mat} and \cite{Lan}.
\begin{lem}\label{decreasing}
Let $0<z\leq \nu$. Then
$$J_{\nu}(z)>0\quad\text{and}\quad\frac{\partial J_{\nu}(z)}{\partial\nu}<0$$
whereas
$$Y_{\nu}(z)<0\quad\text{and}\quad\frac{\partial Y_{\nu}(z)}{\partial\nu}<0.$$

Moreover, as $\xi$ goes to $+\infty$,
\begin{equation}\label{threshold1}
J_{\xi}(\xi)=\frac{\Gamma(1/3)}{3^{1/6}2^{2/3}\pi}\xi^{-1/3}+O(\xi^{-4/3})
\end{equation}
and
\begin{equation}\label{threshold2}
Y_{\xi}(\xi)=-\frac{3^{1/3}\Gamma(1/3)}{2^{2/3}\pi}\xi^{-1/3}+O(\xi^{-4/3}).
\end{equation}
\end{lem}

Using these results we immediately infer the following corollary. The first part is a consequence of Proposition~\ref{firstprop} and the successive computations, the second follows from Lemma~\ref{decreasing}.

\begin{cor}\label{corH=Y}
Let us fix $z_2\geq 1$ and $\varepsilon$, $0<\varepsilon\leq 1/\rme$.
Then there exists a positive constant $\tilde{\nu}_3=\tilde{\nu}_3(\varepsilon)$ such that 
for any $0<z\leq z_2$, and
for any $\nu\geq \max\{\nu_0(\varepsilon/2,z_2),\tilde{\nu}_3(\varepsilon)\}$, we have
\begin{equation}
 |J_{\nu}(z)|\leq \varepsilon |Y_{\nu}(z)|
\end{equation}
hence
$$(1-\varepsilon)|Y_{\nu}(z)|\leq|H^{(1)}_{\nu}(z)|\leq(1+\varepsilon)|Y_{\nu}(z)|.$$

Furthermore, there exist absolute constants $c_1$, $0<c_1<1$, and $\xi_0>0$ such that for any
$\xi\geq \xi_0$ we have
$$ |J_{\xi}(\xi)|\leq c_1 |Y_{\xi}(\xi)|$$
therefore for any $\nu\geq z\geq \xi_0$ we have
$$ |J_{\nu}(z)|\leq c_1 |Y_{\nu}(z)|$$
and hence
\begin{equation}
(1-c_1)|Y_{\nu}(z)|\leq|H^{(1)}_{\nu}(z)|\leq(1+c_1)|Y_{\nu}(z)|.
\end{equation}
\end{cor}

We now consider different regimes. First, we take $\nu$ much larger than $z$, then we consider the case in which $z$ is much larger than $\nu$. Finally we shall deal with the intermediate regime when $\nu$ and $z$ are comparable and large.

We begin with the following two results, due to Matviyenko, proved in Theorem~3.1 and Theorem~3.3 in \cite{Mat}, respectively.

\begin{teo}\label{nu<z}
Let $0\leq\nu<z$. Then
\begin{equation}\label{mat1}
H^{(1)}_{\nu}(z)=\left(\frac{2}{\pi}\right)^{1/2}\frac{1}{(z^2-\nu^2)^{1/4}}\exp(\rmi\eta_1)\left[1+R_1(z,\nu)\right]
\end{equation}
where
$$\eta_1=(z^2-\nu^2)^{1/2}-\nu\arccos(\nu/z)-\pi/4$$
and
$$|R_1(z,\nu)|\leq\exp(\tilde{g}_1)\tilde{g}_1
$$
with
$$g_1=\frac{z-\nu}{z^{1/3}}\quad\text{and}\quad\tilde{g}_1=\frac{2}{3g_1^{3/2}}.$$
\end{teo}

\begin{teo}\label{nu>z}
Let $0<z<\nu$. Then
\begin{equation}\label{mat2}
Y_{\nu}(z)=-\left(\frac{2}{\pi}\right)^{1/2}\frac{1}{(\nu^2-z^2)^{1/4}}\exp(\eta_2)\left[1+R_2(z,\nu)\right]
\end{equation}
where
\begin{equation}\label{eta2def}
\eta_2=\nu\log\left(\frac{\nu}{z}+\left(\left(\frac{\nu}{z}\right)^2-1\right)^{1/2}\right)-(\nu^2-z^2)^{1/2}
\end{equation}
and
$$|R_2(z,\nu)|\leq\exp(\tilde{g}_2)\tilde{g}_2
$$
with
$$g_2=\frac{\nu-z}{\nu^{1/3}}\quad\text{and}\quad\tilde{g}_2=\frac{2}{3g_2^{3/2}}.$$
\end{teo}

We have the following corollaries. By Theorem~\ref{nu<z} we can generalize \eqref{argument}
by extending it uniformly on suitable intervals in the following way.
\begin{cor}
Fixed $C>1$, for any $\nu\geq 0$ and for any positive $z$ such that $z\geq C\nu$ we have
$$\tilde{g}_1\leq \frac{2}{3}\left(\frac{C}{C-1}\right)^{3/2}z^{-1}\quad\text{and}\quad
1\leq \left(1-\frac{\nu^2}{z^2}\right)^{-1/4}\leq \left(\frac{C^2}{C^2-1}\right)^{1/4}.$$
Therefore,
\begin{multline}\label{mat1bis}
\left[1-\frac{2}{3}\left(\frac{C}{C-1}\right)^{3/2}\exp\left(\frac{2}{3}\left(\frac{C}{C-1}\right)^{3/2}z^{-1}\right)z^{-1}\right]\left(\frac{2}{\pi z}\right)^{1/2}
\leq |H^{(1)}_{\nu}(z)| \leq\\
 \left(\frac{2}{\pi z}\right)^{1/2}\left(\frac{C^2}{C^2-1}\right)^{1/4}
\left[1+\frac{2}{3}\left(\frac{C}{C-1}\right)^{3/2}\exp\left(\frac{2}{3}\left(\frac{C}{C-1}\right)^{3/2}z^{-1}\right)z^{-1}\right].
\end{multline}

Furthermore, there exists an absolute constant $C_1>0$ such that for any $\nu\geq 0$ and for any positive $z$ satisfying $z-\nu\geq C_1 z^{1/3}$ then
\begin{equation}\label{mat1tris}
\frac{1}{2}
\left(\frac{2}{\pi}\right)^{1/2}\frac{1}{(z^2-\nu^2)^{1/4}}
\leq |H^{(1)}_{\nu}(z)| \leq
\frac{3}{2}
\left(\frac{2}{\pi}\right)^{1/2}\frac{1}{(z^2-\nu^2)^{1/4}}.
\end{equation}

Finally, there exists a constant $C_2$, depending on $C$ only, such that
for any $\nu\geq 0$ and for any positive $z$ such that $z \geq C\nu$
we have
\begin{equation}\label{argumentgen}
H^{(1)}_{\nu}(z)=\left(\frac{2}{\pi z}\right)^{1/2}\rme^{\rmi(z-\nu\pi/2-\pi/4)}\left[1+\tilde{R}_1\right]
\end{equation}
where
$$|\tilde{R}_1|\leq C_2\left(\frac{\nu^2}{z^2}+\min\left\{2,\frac{\nu^2}{z}\right\}+\exp\left(\frac{2}{3}\left(\frac{C}{C-1}\right)^{3/2}z^{-1}\right)z^{-1}\right).$$
\end{cor}
 
By Theorem~\ref{nu>z} we can generalize \eqref{order}
by extending it uniformly on suitable intervals in the following way.
In fact the following corollary holds and it may be coupled with Corollary~\ref{corH=Y}.

\begin{cor}
Fixed $C>1$, for any $\nu$ such that $\nu\geq Cz>0$ we have
$$\tilde{g}_2\leq \frac{2}{3}\left(\frac{C}{C-1}\right)^{3/2}\nu^{-1}.$$

Thus for any $\nu$ such that $\nu\geq Cz>0$ and $\nu\geq Cz^2$ we have
\begin{equation}\label{ordergen}
H^{(1)}_{\nu}(z)=-\rmi
\left(\frac{2}{\pi\nu}\right)^{1/2}\left(\frac{\rme z}{2\nu}\right)^{-\nu}\left[1+\tilde{R}_2\right]
\end{equation}
where
$$|\tilde{R}_2|\leq \frac{|J_{\nu}(z)|}{|Y_{\nu}(z)|}+C_1\left(1+\frac{|J_{\nu}(z)|}{|Y_{\nu}(z)|}\right)
\left(\frac{z^2}{\nu^2}+\frac{z^2}{\nu}
+\exp\left(\frac{2}{3}\left(\frac{C}{C-1}\right)^{3/2}\nu^{-1}\right)\nu^{-1}\right)$$
and $C_1$ is an absolute constant depending on $C$ only.
\end{cor}

A further important corollary is the following. Here we make use of the continuity properties of the Hankel functions with respect to both the argument and order and the fact that $|H^{(1)}_{\nu}(z)|>0$ for any $\nu\geq 0$ and any $z>0$.
\begin{cor}\label{unifbounded}
Let us fix $0<z_1<z_2$. Then there exists a constant $C\geq 1$, depending on $z_1$ and $z_2$ only, such that for any $z$, $z_1\leq z\leq z_2$, we have
\begin{equation}\label{uniformkbounded0}
C^{-1}\leq|H^{(1)}_0(z)|\leq C
\end{equation}
and for any $\nu\geq 1/2$ we have 
\begin{equation}\label{uniformkboundedl}
C^{-1}\left(\frac{2}{\pi\nu}\right)^{1/2}\left(\frac{\rme z_2}{2\nu}\right)^{-\nu}
\leq C^{-1}\left(\frac{2}{\pi\nu}\right)^{1/2}\left(\frac{\rme z}{2\nu}\right)^{-\nu}
\leq |H^{(1)}_{\nu}(z)|
\end{equation}
and
\begin{equation}\label{uniformkboundedu}
|H^{(1)}_{\nu}(z)|
\leq 
C
\left(\frac{2}{\pi\nu}\right)^{1/2}\left(\frac{\rme z}{2\nu}\right)^{-\nu}
\leq C
\left(\frac{2}{\pi\nu}\right)^{1/2}\left(\frac{\rme z_1}{2\nu}\right)^{-\nu}.
\end{equation}
\end{cor}

This is the estimate that allows us to use Isakov's argument in \cite{Isak92} and prove stability estimates for the determination of the near-fields from far-field measurements, when the wavenumber $k$ belongs to a fixed interval $[k_1,k_2]$, with $0<k_1<k_2$. We shall treat this case at the beginning of Section~\ref{stabsec}, in Theorem~\ref{Isakov}.

We now begin to investigate the more difficult case of the asymptotic behavior of Hankel functions when $z$ and $\nu$ are both large.

We have the following theorem, which is the main result of this section.

\begin{teo}\label{regimes}
There exist positive constants $z_0\geq 1$, $C_0$ and $A_0\geq 1$ such that if $z\geq z_0$ then the following asymptotic behavior of the Hankel function holds for $z>0$ and $\nu\geq 0$.

If
$$\nu>0\text{ and }\frac{z-\nu}{\nu^{1/3}}\geq C_0\quad\text{or}\quad\nu=0$$ 
then
\begin{equation}\label{firstregimeest}
A_0^{-1}\frac{1}{(z^2-\nu^2)^{1/4}}
 \leq|H^{(1)}_{\nu}(z)|\leq A_0\frac{1}{(z^2-\nu^2)^{1/4}}.
\end{equation}

If $\nu>0$ and
$$\frac{|z-\nu|}{\nu^{1/3}}\leq C_0$$ 
then
\begin{equation}\label{secondregimeest}
A_0^{-1}\nu^{-1/3}\leq|H^{(1)}_{\nu}(z)|\leq A_0\nu^{-1/3}.
\end{equation}

If
$$\frac{\nu-z}{\nu^{1/3}}\geq C_0$$ 
then
\begin{equation}\label{thirdregimeest}
A_0^{-1}\exp(\eta_2)\frac{1}{(\nu^2-z^2)^{1/4}} \leq|H^{(1)}_{\nu}(z)|\leq A_0
\exp(\eta_2)\frac{1}{(\nu^2-z^2)^{1/4}}
\end{equation}
where $\eta_2$ is as in \eqref{eta2def}.
\end{teo}

\proof{.} Let us begin with the third regime, when $\nu$ is much greater than $z$.
Let us assume that $\nu> z\geq z_0=\xi_0$. Then \eqref{thirdregimeest} follows immediately by
Corollary~\ref{corH=Y} and Theorem~\ref{nu>z}.

For what concerns the first estimate, \eqref{firstregimeest}, let us notice that it is trivial for $\nu=0$ and $z\geq z_0$. Therefore, without loss of generality, in what follows we shall assume $\nu>0$.
By Theorem~\ref{nu<z},
there exist $C_0$ and $A_0$ such that \eqref{firstregimeest} holds provided
$z-\nu\geq C_0z^{1/3}>0$, with no assumption that $z$ should be greater than a constant.
Clearly $z-\nu\geq C_0z^{1/3}>0$ implies that $z-\nu\geq C_0\nu^{1/3}>0$ but viceversa does not hold. The case $\nu+C_0\nu^{1/3}\leq z <\nu+C_0z^{1/3}$ will follow from the analysis of the intermediate regime, when $z$ and $\nu$ have approximately the same value, which we shall now deal with.

We begin with the following remark. Fixed $C_0>0$, there exists $z_0>0$ such that if $z\geq z_0$, then $z-C_0z^{1/3}\geq z/2$. Therefore, if $z\geq z_0$ and
$0<z-\nu< C_0z^{1/3}$, we also have $\nu<z<\nu+C_02^{1/3}\nu^{1/3}.$
Therefore, the open case is when $z\geq z_0$ and 
$$\nu-C_1\nu^{1/3}\leq z\leq \nu+C_1\nu^{1/3}$$
where $C_1=C_02^{1/3}$. Let us remark here that the intermediate estimate \eqref{secondregimeest} for $\nu=z$ follows from Lemma~\ref{decreasing}, therefore, without loss of generality, we shall assume either $z<\nu$ or $z>\nu$, that is
$$\nu-C_1\nu^{1/3}\leq z<\nu\quad\text{or}\quad\nu<z\leq \nu+C_1\nu^{1/3}.$$

The analysis for these regimes is rather difficult and has been carried out by Langer in \cite{Lan}. We shall here recall this asymptotic analysis.
We need to introduce the following notation
\begin{equation}\label{x>0not}
\begin{array}{ll}
\begin{array}{l}\rme^x=\sec(\beta)\\
\phi(x)=\tan(\beta)\\
\xi=\nu(\tan(\beta)-\beta)\\
\psi(x)=\displaystyle{\frac{(\tan(\beta)-\beta)^{1/6}}{(\tan(\beta))^{1/2}}}
\end{array} &\text{if }x>0,
\end{array}
\end{equation}
and also
\begin{equation}\label{x<0not}
\begin{array}{ll}
\begin{array}{l}\rme^x=\sech(\alpha)\\
\phi(x)=\rmi\tanh(\beta)\\
\xi=\rmi^3\nu(\alpha-\tanh(\alpha))\\
\psi(x)=\displaystyle{\frac{(\alpha-\tanh(\alpha))^{1/6}}{(\tanh(\alpha))^{1/2}}}
\end{array} &\text{if }x< 0.
\end{array}
\end{equation}

Then Langer in \cite{Lan} showed that, fixed a positive constant $M$, there exists a constant $E$ depending on $M$ only such that
if $0<|\xi|\leq M$ then we have
\begin{equation}\label{langerest}
H^{(1)}_{\nu}(\nu\rme^x)=\frac{2\psi(x)\xi^{1/3}}{3^{1/2}\nu^{1/3}}\left[\rme^{-\pi\rmi/3}J_{-1/3}(\xi)+\rme^{\pi\rmi/3}J_{1/3}(\xi)
\right]+\frac{R(x,\nu)}{\nu^{4/3}}
\end{equation}
with
$$|R(x,\nu)|\leq E.$$

First we need to investigate the term $\left[\rme^{-\pi\rmi/3}J_{-1/3}(\xi)+\rme^{\pi\rmi/3}J_{1/3}(\xi)
\right]$. By using the classical asymptotics as $z\to 0$, see for instance \cite[page~44]{Wat}, we have that, for any complex $z\neq 0$,
$$J_{\pm 1/3}(z)=\frac{(z/2)^{\pm1/3}}{\Gamma(1\pm1/3)}(1+R_{\pm})$$
where
$$|R_{\pm}|<\exp\left(\frac{|z|^2/4}{1\pm 1/3}\right)-1.$$
Since for $z\neq 0$ we have
$$J_{1/3}(z)J'_{-1/3}(z)-J_{-1/3}(z)J'_{1/3}(z)=-2\frac{\sin(\pi/3)}{\pi z},$$
we deduce that, for any $|\xi|>0$ defined as before, we have
$$\left|\rme^{-\pi\rmi/3}J_{-1/3}(\xi)+\rme^{\pi\rmi/3}J_{1/3}(\xi)
\right|>0,$$
hence there exist positive constants $B_0<B_1$, depending on $M$ only, such that,
for any $0<|\xi|\leq M$ defined as before, we have
\begin{equation}\label{stima}
B_0|\xi|^{-1/3}\leq \left|\rme^{-\pi\rmi/3}J_{-1/3}(\xi)+\rme^{\pi\rmi/3}J_{1/3}(\xi)
\right|\leq B_1|\xi|^{-1/3}.
\end{equation}
Then we have that
$$\frac{2B_0}{\nu^{1/3}}\frac{|\psi(x)|}{3^{1/2}}\leq\left|\frac{2\psi(x)\xi^{1/3}}{3^{1/2}\nu^{1/3}}\left[\rme^{-\pi\rmi/3}J_{-1/3}(\xi)+\rme^{\pi\rmi/3}J_{1/3}(\xi)
\right]\right|\leq \frac{2B_1}{\nu^{1/3}}\frac{|\psi(x)|}{3^{1/2}}.
$$

It remains to analyze the term $|\psi(x)|/3^{1/2}$. We need to consider separately the case $x>0$ and $x<0$. For $x>0$, we have $0<\beta<\pi/2$ and
$$\frac{|\psi(x)|}{3^{1/2}}=\displaystyle{\frac{(\tan(\beta)-\beta)^{1/6}}{(3\tan(\beta))^{1/2}}}.$$
There exists an absolute constant $\beta_0$, $0<\beta_0<\pi/2$ such that for any $\beta$, $0\leq\beta\leq \beta_0$ we have
$$\beta+\beta^3/3\leq\tan(\beta)\leq \beta+2\beta^3/3\leq (3/2)\beta\quad\text{and}\quad 1+\beta^2/4\leq\sec(\beta)\leq 1+3\beta^2/4.$$

Let $\nu<z<\nu+C_1\nu^{1/3}$, then $\sec(\beta)=z/\nu=1+a\nu^{-2/3}$, for some $a$, $0<a\leq C_1$. There exists $\tilde{z}_0>0$, depending on $C_1$ only, such that if $z\geq \tilde{z}_0$ then $\beta\leq\beta_0$. Hence
$$(2/\sqrt{3})\sqrt{a}\nu^{-1/3}\leq\beta\leq 2\sqrt{a}\nu^{-1/3}$$
that is
$$\frac{8}{9\sqrt{3}}a^{3/2}\leq \nu(\tan(\beta)-\beta)=\xi\leq \frac{16}{3}a^{3/2}\leq \frac{16}{3}C_1^{3/2}=M,$$
with $M$ thus depending on $C_1$ only. An easy computation shows that there exist absolute positive constants $\tilde{B}_0<\tilde{B}_1$ such that
$$\tilde{B}_0\leq \displaystyle{\frac{(\tan(\beta)-\beta)^{1/6}}{(3\tan(\beta))^{1/2}}}\leq\tilde{B}_1.$$
We may conclude that there exist positive constants $\tilde{z}_0$ and $\tilde{A}_0< \tilde{A}_1$, depending on $C_1$ only, such that for any $z\geq \tilde{z}_0$ and
$\nu<z<\nu+C_1\nu^{1/3}$ we have
$$\frac{\tilde{A}_0}{\nu^{1/3}}\leq\left|\frac{2\psi(x)\xi^{1/3}}{3^{1/2}\nu^{1/3}}\left[\rme^{-\pi\rmi/3}J_{-1/3}(\xi)+\rme^{\pi\rmi/3}J_{1/3}(\xi)
\right]\right|\leq \frac{\tilde{A}_1}{\nu^{1/3}}.
$$
Since the constant $E$ depends on $M$, thus on $C_1$ only, we can find $z_0\geq\tilde{z}_0$ such that for any $z\geq z_0$ and
$\nu<z<\nu+C_1\nu^{1/3}$ we have
$$E\nu^{-1}\leq (1/2)\tilde{A}_0,$$ 
thus \eqref{secondregimeest} is proved with $A_0=\tilde{A}_0/2$ and $A_1=2\tilde{A}_1$.
A completely analogous argument holds for the case $x<0$, hence \eqref{secondregimeest} is fully proved.

It remains to consider the case when $z\geq z_0$ and
$\nu+C_0\nu^{1/3}\leq z <\nu+C_0z^{1/3}\leq \nu+C_1\nu^{1/3}$.
We have just obtained that
$A_0^{-1}\nu^{-1/3}\leq|H^{(1)}_{\nu}(z)|\leq A_0\nu^{-1/3}$.
On the other hand, provided $z=\nu+a\nu^{1/3}$ with $C_0\leq a\leq C_1$, 
we have that
$$\frac{1}{(z^2-\nu^2)^{1/4}}=\frac{1}{(2a\nu^{4/3}+a^2\nu^{2/3})^{1/4}}=
\nu^{-1/3}\frac{1}{(2a+a^2\nu^{-2/3})^{1/4}}
$$
that is 
$$
\frac{1}{(2C_1+C_1^2\nu^{-2/3})^{1/4}}
\nu^{-1/3}\leq 
\frac{1}{(z^2-\nu^2)^{1/4}}\leq \frac{1}{(2C_0)^{1/4}}\nu^{-1/3}.$$
Provided $z_0$ is big enough, depending on $C_0$ and $C_1$ only, then if $z\geq z_0$ we infer that $\nu\geq 1$, thus
$$
\frac{1}{(2C_1+C_1^2)^{1/4}}
\nu^{-1/3}\leq 
\frac{1}{(z^2-\nu^2)^{1/4}}\leq \frac{1}{(2C_0)^{1/4}}\nu^{-1/3}$$
and the proof is concluded.\cvd

\section{Stability estimates: from far-field to near-field}\label{stabsec}

Throughout this section we shall assume that $u^s$ is a radiating solution of the Helmholtz equation $\Delta u^s+k^2u^s=0$, for some wavenumber $k>0$, defined in $\mathbb{R}^N\backslash \overline{B_R}$, for some fixed $R>0$. We call $u^s_{\infty}$ the corresponding far-field pattern.

Throughout this section
we shall also fix $R_0>R$ and $B_0$ and $B_1$, with $1<B_0<B_1$. We call $\tilde{b}_0=1/B_0$ and we assume that for some $\varepsilon>0$ and $M>0$ we have
\begin{equation}\label{farfielderrormeas}
\|u^s_{\infty}\|_{L^2(\mathbb{S}^{N-1})}=\varepsilon,\quad\|u^s\|_{L^2(\partial B_{R_0})}=M.\end{equation}

Our aim is to estimate, in terms of $\|u^s_{\infty}\|_{L^2(\mathbb{S}^{N-1})}$, $\|u^s\|_{L^2(\partial B_r)}$ for some $r>R$. This kind of estimate is usually referred to as the stability for the determination of the near-field from the far-field. Notice the Rellich Lemma provides the corresponding uniqueness. Such an issue has been solved by Isakov, \cite{Isak92}, see also \cite{Bus}.

Let us recall there exists a sequence $b_{j}$, $j=0,1,\ldots,$ of nonnegative numbers such that
\begin{equation}\label{ff}
\|u^s_{\infty}\|^2_{L^2(\mathbb{S}^{N-1})}=\sum_jb_j^2
\end{equation}
and for any $r>R$
\begin{equation}\label{nf}
\|u^s\|^2_{L^2(\partial B_r)}= \frac{\pi}{2}
\sum_jb_j^2
kr
\left|H^{(1)}_{j+(N-2)/2}(kr)\right|^2.
\end{equation}

Let us notice that for any $j\geq 0$,
$$b_j^2=\sum_{i:\ \gamma(v_i)=j}|\tilde{b}_i|^2,$$
if we use the notation of Section~\ref{farfield} as in \eqref{decompscat3} and \eqref{farfielderror}.

For any $k>0$ and any $r$, $B_0R_0\leq r\leq B_1R_0$, we have,
for any integer $j_0\geq 0$ to be decided later,
\begin{equation}\label{nf2}
\|u^s\|^2_{L^2(\partial B_r)}\leq \frac{\pi}{2}kr
\sum_{j=0}^{j_0}b_j^2
\left|H^{(1)}_{j+(N-2)/2}(kr)\right|^2+\frac{\pi}{2}kr
\sum_{j=j_0+1}^{+\infty}b_j^2
\left|H^{(1)}_{j+(N-2)/2}(kr)\right|^2
\end{equation}
which we estimate as follows
\begin{multline}\label{nf3}
\|u^s\|^2_{L^2(\partial B_r)}\leq \frac{\pi}{2}kr
\max_{j\in\{0,1,\ldots,j_0\}}\left|H^{(1)}_{j+(N-2)/2}(kr)\right|^2\|u^s_{\infty}\|^2_{L^2(\mathbb{S}^{N-1})}
+\\
\frac{r}{R_0}\frac{\pi}{2}kR_0
\sum_{j=j_0+1}^{+\infty}b_j^2\left|H^{(1)}_{j+(N-2)/2}(kR_0)\right|^2
\frac{\left|H^{(1)}_{j+(N-2)/2}(kr)\right|^2}{\left|H^{(1)}_{j+(N-2)/2}(kR_0)\right|^2}
\end{multline}
that is
\begin{multline}\label{nf4}
\|u^s\|^2_{L^2(\partial B_r)}\leq \frac{\pi}{2}kr
\max_{j\in\{0,1,\ldots,j_0\}}\left|H^{(1)}_{j+(N-2)/2}(kr)\right|^2\|u^s_{\infty}\|^2_{L^2(\mathbb{S}^{N-1})}
+\\
\frac{r}{R_0}\left(\sup_{j>j_0}
\frac{\left|H^{(1)}_{j+(N-2)/2}(kr)\right|^2}{\left|H^{(1)}_{j+(N-2)/2}(kR_0)\right|^2}
\right)\|u^s\|^2_{L^2(\partial B_{R_0})}.
\end{multline}

We recall here the stability result by Isakov, \cite{Isak92}, which we slightly generalize to any dimension $N\geq 2$ and for the wavenumber $k$ varying in a compact interval contained in $(0,+\infty)$. In fact, let us fix two positive constants $k_1$, $k_2$ with $0<k_1<k_2$.
Let us fix $k\in [k_1,k_2]$. We call $z_1=k_1R_0$ and $z_2=k_2B_1R_0$.
By Corollary~\ref{unifbounded},
for any $\nu$ such that $\nu=0$ or $1/2\leq\nu\leq \rme z/2$, we have that there exists a constant $C\geq 1$, depending on $z_1$ and $z_2$ only, such that
\begin{equation}
|H^{(1)}_{\nu}(z)|
\leq 
C\frac{2}{\sqrt{\rme\pi z}}.
\end{equation}
Otherwise, if $1/2\leq\nu$ and $\rme z/2\leq \nu$, we have
\begin{equation}
|H^{(1)}_{\nu}(z)|
\leq 
C\frac{2}{\sqrt{\rme\pi z}}\left(\frac{2\nu}{\rme z}\right)^{\nu-1/2}.
\end{equation}
Furthermore, for any $\nu\geq 1/2$ we have
\begin{equation}
\frac{\left|H^{(1)}_{\nu}(kr)\right|}{\left|H^{(1)}_{\nu}(kR_0)\right|}\leq
\frac{C\frac{2}{\sqrt{\rme\pi kr}}\left(\frac{2\nu}{\rme kr}\right)^{\nu-1/2}}{C^{-1}\frac{2}{\sqrt{\rme\pi kR_0}}\left(\frac{2\nu}{\rme kR_0}\right)^{\nu-1/2}}\leq C^2\left(\frac{R_0}{r}\right)^{\nu}\leq
C^2\left(1/B_0\right)^{\nu}.
 \end{equation}

We conclude that, provided $j_0\geq \hat{j}_0\geq 1$, with $\hat{\nu}_0=\hat{j}_0+(N-2)/2\geq\max\{N/2,\rme^2 z_2/2\}$,
and setting $\nu_0=j_0+(N-2)/2$,
\begin{equation}\label{nf6}
\|u^s\|^2_{L^2(\partial B_r)}\leq \frac{2C^2}{\rme}
\left(\frac{2\nu_0}{\rme kr}\right)^{2\nu_0-1}
\|u^s_{\infty}\|^2_{L^2(\mathbb{S}^{N-1})}
+C^4\left(\frac{R_0}{r}\right)^{2\nu_0-1}\|u^s\|^2_{L^2(\partial B_{R_0})},
\end{equation}
that is, setting $\tilde{A}=\max\{2C^2/\rme,C^4\}$ and $\tilde{a}_1$, $0<\tilde{a}_1<1$, such that
$\rme^{\tilde{a}_1}=\rme/2$, and $\tilde{b}_0=1/B_0$, 
\begin{equation}\label{nf8}
\|u^s\|^2_{L^2(\partial B_r)}\leq
\tilde{A}\left[
\left(\frac{2\nu_0-1}{\rme^{\tilde{a}_1} kr}\right)^{2\nu_0-1}
\|u^s_{\infty}\|^2_{L^2(\mathbb{S}^{N-1})}
+\left(\tilde{b}_0\right)^{2\nu_0-1}\|u^s\|^2_{L^2(\partial B_{R_0})}\right].
\end{equation}

Let $2n\geq 2\hat{\nu}_0+1$
be such that
\begin{equation}\label{choicen0}
\left(\frac{2n}{\rme^{\tilde{a}_1}kr}\right)^{2n}\varepsilon^2= \tilde{b}_0^{2n}M^2.
\end{equation}
Then there exists an integer $j_0\in \mathbb{N}$
such that $j_0\geq \hat{j}_0$ and
$2\nu_0-1\leq 2n< 2(j_0+1+(N-2)/2)-1=2\nu_0+1$.
Then 
$$\left(\frac{2\nu_0-1}{\rme^{\tilde{a}_1}kr}\right)^{2\nu_0-1}\leq \left(\frac{2n}{\rme^{\tilde{a}_1}kr}\right)^{2n}$$
and
$$\tilde{b}_0^{2\nu_0-1}=\tilde{b}_0^{-2}\tilde{b}_0^{2\nu_0+1}\leq \tilde{b}_0^{-2}\tilde{b}_0^{2n}.$$
Hence, if \eqref{choicen0} holds we have 
\begin{equation}
\|u^s\|^2_{L^2(\partial B_r)}\leq 
2\tilde{A}\tilde{b}_0^{-2}
\left(\frac{2n}{\rme^{\tilde{a}_1}kr}\right)^{2n}\varepsilon^2
=
2\tilde{A}\tilde{b}_0^{-2}
\tilde{b}_0^{2n}M^2.
\end{equation}

Let us investigate \eqref{choicen0}.
We have
$\log(\varepsilon)+n\log\left(\frac{2n}{\rme^{\tilde{a}_1}kr}\right)=\log(M)+n\log(\tilde{b}_0)$, that is
$$\log(M/\varepsilon)=n\log\left(\frac{2n}{\tilde{b}_0\rme^{\tilde{a}_1}kr}\right).$$
Let us call
$$\tilde{n}=\frac{2n}{\tilde{b}_0\rme^{\tilde{a}_1}kr}
,$$
then
$$\frac{2}{\tilde{b}_0\rme^{\tilde{a}_1}kr}\log(M/\varepsilon)=\tilde{n}\log(\tilde{n}).$$
We notice that if $2n\geq 2\hat{\nu}_0+1$, then $\tilde{n}\geq \rme$.
On the other hand, for any $\alpha>0$, there exists $C_{\alpha}\geq 1$ such that for any $\tilde{n}\geq \rme$ we have
$$\tilde{n}\leq \tilde{n}\log(\tilde{n})\leq C_{\alpha}\tilde{n}^{1+\alpha}.$$
Notice that we can choose $C_{\alpha}=1$ for any $\alpha\geq \alpha_0$ where
$\alpha_0$
satisfies $0<\alpha_0<1$ and 
$$(1/\alpha_0)^{1/\alpha_0}=\rme.$$
Therefore,
$$\tilde{n}\leq
\frac{2}{\tilde{b}_0\rme^{\tilde{a}_1}kr}\log(M/\varepsilon)
\leq C_{\alpha}\tilde{n}^{1+\alpha},
$$
that is
$$\tilde{n}\geq \left(\frac{2}{\tilde{b}_0\rme^{\tilde{a}_1}kr}\frac{\log(M/\varepsilon)}{C_{\alpha}}\right)^{1/(1+\alpha)}.$$

Therefore, provided
$$
\left(\frac{2}{\tilde{b}_0\rme^{\tilde{a}_1}kr}\frac{\log(M/\varepsilon)}{C_{\alpha}}\right)^{1/(1+\alpha)}\geq \frac{2\hat{\nu}_0+1}{\tilde{b}_0\rme^{\tilde{a}_1}kr},
$$
that is
\begin{equation}\label{cond}
\log(M/\varepsilon)\geq C_{\alpha}
\left(\frac{2}{\tilde{b}_0\rme^{\tilde{a}_1}kr}\right)^{\alpha}
(\hat{\nu}_0+1/2)^{1+\alpha},
\end{equation}
we have that there exists a solution $n$ to \eqref{choicen0} such that $2n\geq 2\hat{\nu}_0+1$,
hence
$$
\|u^s\|_{L^2(\partial B_r)}\leq 
\sqrt{2\tilde{A}}\tilde{b}_0^{-1}
M\exp (-\log(B_0)n),$$
that is
\begin{equation}\label{est}
\|u^s\|_{L^2(\partial B_r)}\leq
\sqrt{2\tilde{A}}\tilde{b}_0^{-1}
M\exp\left(-\log(B_0)
\left(\left(\frac{\tilde{b}_0\rme^{\tilde{a}_1}kr}{2}\right)^{\alpha}\frac{\log(M/\varepsilon)}{C_{\alpha}}\right)^{1/(1+\alpha)}\right).
\end{equation}

We have just proved the following result.

\begin{teo}\label{Isakov}
Let $N\geq 2$. Let us fix constants $k_1$, $k_2$ with $0<k_1<k_2$. Under the previous notation and assumptions, the following result holds.

Let $z_1=k_1R_0$ and $z_2=k_2B_1R_0$ and let $\hat{\nu}_0=\max\{N/2,\rme^2z_2/2\}$. Let $\tilde{a}_1=\log(\rme/2)$.
Fix $\alpha>0$. Then there exist a constant $C_{\alpha}\geq 1$, depending on $\alpha$ only \textnormal{(}with $C_{\alpha}=1$ for any $\alpha\geq \alpha_0$\textnormal{)}, and a constant $\tilde{A}$, depending on $z_1$ and $z_2$ only, such that
for any $k\in [k_1,k_2]$ and any $r$, $B_0R_0\leq r\leq B_1R_0$,
if
\begin{equation}\label{cond2}
\log(M/\varepsilon)\geq C_{\alpha}
\left(\frac{2}{\rme^{\tilde{a}_1}k_1R_0}\right)^{\alpha}
(\hat{\nu}_0+1/2)^{1+\alpha},
\end{equation}
then
\begin{equation}\label{est2}
\|u^s\|_{L^2(\partial B_r)}\leq 
\sqrt{2\tilde{A}}\tilde{b}_0^{-1}
M\exp\left(-\log(B_0)
\left(\left(\frac{\rme^{\tilde{a}_1}k_1R_0}{2}\right)^{\alpha}\frac{\log(M/\varepsilon)}{C_{\alpha}}\right)^{1/(1+\alpha)}\right).
\end{equation}
Obviously, we can replace \eqref{cond2} and \eqref{est2} with \eqref{cond} and \eqref{est}.
\end{teo}

Now we investigate how the estimate changes in the high frequencies regime. We shall keep the notation set at the beginning of the section.
Let $z_0$, $C_0$ and $A_0$ be as in 
Theorem~\ref{regimes}.
We observe that there exists a constant $\tilde{C}_0\geq 2$, depending on $C_0$ only, such that for any $\nu\geq 1/2$ we have
$\nu+C_0\nu^{1/3}\leq \tilde{C}_0\nu$. 
Without loss of generality, 
up to taking a possibly greater $z_0$ but still depending only on absolute constants and on $N$, we
assume also that
\begin{equation}\label{extra2}
z_0\geq \tilde{C}_0\max\{2,(N-2)/2\}\quad\text{and}\quad(1-C_0z_0^{-2/3})>0.
\end{equation}
We further assume that 
\begin{equation}\label{extra}
\text{if }z_0\leq \tilde{C}_0\nu\text{ then
 }C_0\nu^{1/3}\leq \nu/2.
\end{equation}
Notice that if $z_0\leq \tilde{C}_0\nu$ we also have $\nu\geq 2$.


Let us fix $b_0$, $0<b_0<1$, depending on $B_0$ only, such that $1/B_0<b_0<1$.
We fix a positive constant $C\geq \rme^2/2$ and let $a=\sqrt{C^2-1}/C$, $0<a<1$.
We assume that $C(1-C_0z_0^{-2/3})\geq 1$, that $\frac{2\rme^{1-a}}{(1+a)B_0}\leq b_0$ and that $a\geq 2\log(4/3)$.
Obviously $C$ depends on $B_0$, $C_0$ and $z_0$ only.

We fix a constant
$k_0>0$ such that $k_0R_0\geq z_0$. Let us consider $r$, $B_0R_0\leq r\leq B_1R_0$, and $k\geq k_0$. 

Let us notice that if $\nu\geq Ckr\geq Cz_0$ then
$\nu-C_0\nu^{1/3}\geq \nu(1-C_0\nu^{-2/3})\geq \nu(1-C_0z_0^{-2/3})$.
Therefore
$\nu-C_0\nu^{1/3}\geq Ckr(1-C_0z_0^{-2/3})\geq kr$.

Let $\hat{j}_0=\hat{j}_0(kr)\geq 1$ be such that $\hat{\nu}_0=\hat{\nu}_0(kr)=\hat{j}_0(kr)+(N-2)/2\geq\max\{N/2,Ckr\}\geq 2$.
Let $2n\geq 2\hat{\nu}_0(kr)+1$, then
there exists an integer $j_0\in \mathbb{N}$
such that $j_0\geq \hat{j}_0(kr)$ and
$2\nu_0-1\leq 2n< 2(j_0+1+(N-2)/2)-1=2\nu_0+1$, $\nu_0=j_0+(N-2)/2$ as before.

We consider an integer $j_0\geq \hat{j}_0$, to be fixed later, and we use Theorem~\ref{regimes}. If $\nu=0$ or $\nu\geq 1/2$ and such that
$kr-\nu\geq C_0\nu^{1/3}$, we have that, by \eqref{firstregimeest},
$$\left|H^{(1)}_{\nu}(kr)\right|\leq A_0\frac{1}{((kr)^2-\nu^2)^{1/4}}\leq A_0
\max\left\{1/\sqrt{z_0},\frac{1}{(2C_0(1/2)^{4/3}+C_0^2(1/2)^{2/3})^{1/4}}\right\}=\tilde{A}_1$$
since, for $\nu=0$, $1/\sqrt{kr}\leq 1/\sqrt{z_0}$, and for $\nu\geq 1/2$,
$$((kr)^2-\nu^2)\geq (2C_0\nu^{4/3}+C_0^2\nu^{2/3})\geq (2C_0(1/2)^{4/3}+C_0^2(1/2)^{2/3}).$$
Clearly $\tilde{A}_1$ depends on $A_0$, $z_0$ and $C_0$ only.

If $\nu\geq 1/2$ and
$|kr-\nu|\leq C_0\nu^{1/3}$, then, by \eqref{secondregimeest},
$$\left|H^{(1)}_{\nu}(kr)\right|\leq A_0\nu^{-1/3}\leq A_0(1/2)^{-1/3}=\tilde{A}_2,$$
with $\tilde{A}_2$ depending on $A_0$ only.

If $1/2\leq \nu\leq \nu_0$ and
$\nu-kr\geq C_0\nu^{1/3}$, then, by Corollary~\ref{corH=Y} and Lemma~\ref{decreasing} we have that
$$\left|H^{(1)}_{\nu}(kr)\right|\leq \frac{1+c_1}{1-c_1}|H^{(1)}_{\nu_0}(kr)|.$$

For any $\nu\geq \nu_0$ we can use \eqref{thirdregimeest} both for $z=kR_0$ and $z=kr$. Recalling that
$a=\sqrt{C^2-1}/C$, $0<a<1$, we obtain that
\begin{equation}\label{crucial}
A_0^{-1}\frac{1}{\sqrt{\nu}}\left(\frac{(1+a)\nu}{\rme z}\right)^{\nu}
\leq \left|H^{(1)}_{\nu}(z)\right|\leq A_0\frac{1}{\sqrt{a\nu}}
\left(\frac{2\nu}{\rme^{a}z}\right)^{\nu}.
\end{equation}
Hence, if $1/2\leq \nu\leq \nu_0$ and
$\nu-kr\geq C_0\nu^{1/3}$, we have
$$\left|H^{(1)}_{\nu}(kr)\right|\leq \frac{1+c_1}{1-c_1}
A_0\frac{1}{\sqrt{a\nu_0}}
\left(\frac{2\nu_0}{\rme^{a}kr}\right)^{\nu_0}.$$
Since $\nu_0\geq 2$, we have that $2\nu_0\leq (4/3)(2\nu_0-1)$. Let $a_1=a-\log(4/3)$, then
$$\left|H^{(1)}_{\nu}(kr)\right|\leq \tilde{A}_3\left(\frac{\rme^{a_1}}{2\nu_0-1}\right)^{1/2}
\left(\frac{2\nu_0-1}{\rme^{a_1}kr}\right)^{\nu_0},$$
with $\tilde{A}_3$ depending on $A_0$ and $C$ only.


We conclude that for any $j_0\geq \hat{j}_0(kr)$,
setting $A_1=\max\{\tilde{A}_1\sqrt{z_0},\tilde{A}_2\sqrt{z_0},\tilde{A}_3\}$, we have
\begin{equation}
\max_{j\in\{0,1,\ldots,j_0\}}\left|H^{(1)}_{j+(N-2)/2}(kr)\right|^2\leq A^2_1
\left(\frac{\rme^{a_1}}{2\nu_0-1}\right)
\left(\frac{2\nu_0-1}{\rme^{a_1}kr}\right)^{2\nu_0}.
\end{equation}

Now we investigate the integers $j> j_0$. Again by \eqref{crucial}, for any $\nu\geq \nu_0$
\begin{equation}
\frac{\left|H^{(1)}_{\nu}(kr)\right|}{\left|H^{(1)}_{\nu}(kR_0)\right|}\leq
\frac{A_0^2}{\sqrt{a}} \left(\frac{2\rme^{1-a} R_0}{(1+a)r}\right)^{\nu}\leq
\frac{A_0^2}{\sqrt{a}} \left(\frac{2\rme^{1-a}}{(1+a)B_0}\right)^{\nu}\leq \frac{A_0^2}{\sqrt{a}}(b_0)^{\nu}=A_2(b_0)^{\nu-1/2}
\end{equation}
where $A_2=\sqrt{b_0}A_0^2/\sqrt{a}$.
We may conclude that, setting $A=\max\{(\pi/2)A^2_1,A^2_2\}$ and $\tilde{B}=(N/2)C$, if
$$\hat{\nu}_0(kr)=\tilde{B}kr\quad\text{and}\quad \nu_0\geq \hat{\nu}_0(kr),$$
then we have
\begin{equation}
\|u^s\|^2_{L^2(\partial B_r)}\leq 
AB_1
\left[
\left(\frac{2\nu_0-1}{\rme^{a_1}kr}\right)^{2\nu_0-1}
\|u^s_{\infty}\|^2_{L^2(\mathbb{S}^{N-1})}+
\left(b_0\right)^{2\nu_0-1}\|u^s\|^2_{L^2(\partial B_{R_0})}
\right].
\end{equation}
Notice that $A$ depends on $B_0$, $C_0$, $z_0$ and $A_0$ only,
$\tilde{B}$ depends on $B_0$, $C_0$, $z_0$ and $N$ only, $a_1$ depends on $B_0$, $C_0$ and $z_0$ only,
whereas
$0<b_0<1$ depends on $B_0$ only. 

Let $2n\geq 3\hat{\nu}_0(kr)\geq 2\hat{\nu}_0(kr)+1$ be such that
\begin{equation}\label{choicen}
\left(\frac{2n}{\rme^{a_1}kr}\right)^{2n}\varepsilon^2= b_0^{2n}M^2.
\end{equation}
Let $j_0\in \mathbb{N}$ be such that $2\nu_0-1\leq 2n< 2\nu_0+1$.
Then $j_0\geq \hat{j}_0(kr)$ and
$$\left(\frac{2\nu_0-1}{\rme^{a_1}kr}\right)^{2\nu_0-1}\leq \left(\frac{2n}{\rme^{a_1}kr}\right)^{2n}$$
and
$$b_0^{2\nu_0-1}=b_0^{-2}b_0^{2\nu_0+1}\leq b_0^{-2}b_0^{2n}.$$
Hence, if \eqref{choicen} holds we have 
\begin{equation}
\|u^s\|^2_{L^2(\partial B_r)}\leq 
2AB_1b_0^{-2}
\left(\frac{2n}{\rme^{a_1}kr}\right)^{2n}\varepsilon^2
=
2AB_1b_0^{-2}
b_0^{2n}M^2.
\end{equation}

We argue exactly as in the previous case, with $\tilde{A}$, $\tilde{b}_0$ and $\tilde{a}_1$ replaced by $AB_1$, $b_0$ and $a_1$, respectively.
Fixed $\alpha>0$ and $C_{\alpha}$ as before, provided
\begin{equation}\label{cond3}
\log(M/\varepsilon)\geq C_{\alpha}
\left(\frac{2}{b_0\rme^{a_1}kr}\right)^{\alpha}
((3/2)\hat{\nu}_0(kr))^{1+\alpha},
\end{equation}
we have that there exists a solution $n$ to \eqref{choicen} such that $2n\geq 3\hat{\nu}_0(kr)\geq  2\hat{\nu}_0(kr)+1$,
hence
\begin{equation}\label{est3}
\|u^s\|_{L^2(\partial B_r)}\leq 
\sqrt{2AB_1}b_0^{-1}
M\exp\left(-\log(1/b_0)
\left(\left(\frac{b_0\rme^{a_1}kr}{2}\right)^{\alpha}\frac{\log(M/\varepsilon)}{C_{\alpha}}\right)^{1/(1+\alpha)}\right).
\end{equation}

%
%
%

Clearly such an estimate improves as $k$ becomes larger. However, if we call
$$k_1(\varepsilon,r)=\frac{1}{r}\left(\frac{2}{3\tilde{B}}\right)^{1+\alpha}\left(\frac{b_0\rme^{a_1}}{2}\right)^{\alpha}\frac{\log(M/\varepsilon)}{C_{\alpha}},
$$
and we assume that $k_1(\varepsilon,r)\geq k_0$, then we obtain that the estimate
remains valid in the following regime
\begin{equation}\label{regimesk}
k_0\leq k\leq k_1(\varepsilon,r)=\frac{1}{r}
\left(\frac{2}{3\tilde{B}}\right)^{1+\alpha}\left(\frac{b_0\rme^{a_1}}{2}\right)^{\alpha}\frac{\log(M/\varepsilon)}{C_{\alpha}}.
\end{equation}
If we pick the optimal choice of $k$, that is $k=k_1(\varepsilon,r)$, then \eqref{est3} reduces to
\begin{multline*}
\|u^s\|_{L^2(\partial B_r)}\leq 
\sqrt{2AB_1}b_0^{-1}
M\exp\left(-\log(1/b_0)(3\tilde{B}/2)k_1(\varepsilon,r)r\right)=\\
\sqrt{2AB_1}b_0^{-1}
M\exp\left(-\log(1/b_0)
\left(\frac{b_0\rme^{a_1}}{3\tilde{B}}\right)^{\alpha}
\frac{\log(M/\varepsilon)}{C_{\alpha}}\right)
\end{multline*}
that is
\begin{equation}\label{finalestimate2}
\|u^s\|_{L^2(\partial B_r)}\leq \sqrt{2AB_1}b_0^{-1}M^{1-\beta}\varepsilon^{\beta}
\end{equation}
where $\beta=\beta(\alpha)$ is given by
\begin{equation}\label{beta}
\beta(\alpha)=\frac{\log(1/b_0)}{C_{\alpha}}
\left(\frac{b_0\rme^{a_1}}{3\tilde{B}}\right)^{\alpha}
\end{equation}
and, in the particular case $\alpha=1$,
\begin{equation}\label{beta_1}
\beta(1)=\log(1/b_0)
\frac{b_0\rme^{a_1}}{3\tilde{B}}.
\end{equation}
Let us note that, without loss of generality by taking an eventually larger $\tilde{B}$ depending on $\alpha$ as well,  we may assume that $0< \beta(\alpha)\leq 1$. In particular, we have proved the following result, the main one of this section.

\begin{teo}\label{mainteo}
Let $N\geq 2$. We keep the previous notation and assumptions.
Let
$k_0>0$ be such that $k_0R_0\geq z_0$, $z_0$ as in Theorem~\textnormal{\ref{regimes}}
and such that \eqref{extra2} and \eqref{extra} are satisfied. 
Let us fix $b_0$, $0<b_0<1$, depending on $B_0$ only, such that $1/B_0<b_0<1$. Fix $\alpha>0$. 

Then there exist a constant $\tilde{B}$ depending on
$B_0$, $C_0$, $z_0$ and $N$ only, a constant $a_1$, $0<a_1<1$, depending on
$B_0$, $C_0$ and $z_0$ only, a constant $A$ depending on $B_0$, $C_0$, $z_0$ and $A_0$ only, and a constant $C_{\alpha}\geq 1$, depending on $\alpha$ only \textnormal{(}with $C_{\alpha}=1$ for any $\alpha\geq \alpha_0$\textnormal{)}, 
such that the following holds.

Let us assume that
$$k_0\leq k_1(\varepsilon)=\frac{1}{B_1R_0}\left(\frac{2}{3\tilde{B}}\right)^{1+\alpha}\left(\frac{b_0\rme^{a_1}}{2}\right)^{\alpha}\frac{\log(M/\varepsilon)}{C_{\alpha}}.$$
Then for any $k\in [k_0,k_1(\varepsilon)]$ and any $r$, $B_0R_0\leq r\leq B_1R_0$,
we have
\begin{equation}\label{finalest}
\|u^s\|_{L^2(\partial B_r)}\leq \sqrt{2AB_1}b_0^{-1}
M\exp\left(-\log(1/b_0)
\left(\left(\frac{b_0\rme^{a_1}kr}{2}\right)^{\alpha}\frac{\log(M/\varepsilon)}{C_{\alpha}}\right)^{1/(1+\alpha)}\right).
\end{equation}
Furthermore, if $k=k_1(\varepsilon)$ we have
\begin{equation}\label{holder}
\|u^s\|_{L^2(\partial B_r)}
\leq \sqrt{2AB_1}b_0^{-1}
M^{1-\beta_1}\varepsilon^{\beta_1},
\end{equation}
where
$\beta_1=\beta_1(\alpha)$ is given by
\begin{equation}\label{beta1}
\beta_1(\alpha)=\frac{\log(1/b_0)}{C_{\alpha}}
\left(\frac{b_0\rme^{a_1}}{3\tilde{B}}\right)^{\alpha}\left(\frac{r}{B_1R_0}\right)^{\alpha/(1+\alpha)}=
\tilde{C}_{\alpha}
\left(\frac{r}{B_1R_0}\right)^{\alpha/(1+\alpha)}\leq 
\tilde{C}_{\alpha}
\end{equation}
and, in the particular case $\alpha=1$,
\begin{equation}\label{beta1_1}
\beta_1(1)=\log(1/b_0)
\frac{b_0\rme^{a_1}}{3\tilde{B}}\left( \frac{r}{B_1R_0}\right)^{1/2}=
\tilde{C}_1\left( \frac{r}{B_1R_0}\right)^{1/2}
\leq 
\tilde{C}_1.
\end{equation}
Here $\tilde{C}_{\alpha}$ depends on $\alpha$, $B_0$, $C_0$, $z_0$ and $N$ only, and,
without loss of generality by taking an eventually larger $\tilde{B}$ depending on $\alpha$ as well,  we may assume that $0< \beta_1(\alpha)\leq\tilde{C}_{\alpha}
\leq 1$. 
\end{teo}

In the final part of this section, we are interested in understanding what happens if $k>k_1(\varepsilon,r)$. We need to make an additional assumption, namely that
\begin{equation}\label{B_0ass}
B_0\geq \tilde{B}_0=3\rme\tilde{C}_0/2\geq 3\rme.
\end{equation}

We recall that $k_0>0$ is such that $k_0R_0\geq z_0$. We fix $r$, $B_0R_0\leq r\leq B_1R_0$, and $k\geq k_0$. 

We consider all nonnegative integers $j$ such that
$\nu=j+(N-2)/2$ satisfies $kr\geq \tilde{C}_0 \nu$ and we call $\tilde{j}_0=\tilde{j}_0(kr)$ the largest of these integers and $\tilde{\nu}_0=\tilde{\nu}_0(kr)=\tilde{j}_0(kr)+(N-2)/2$. Notice that
$\tilde{j}_0(kr)\geq 0$ by our assumption \eqref{extra2} on $z_0$.
Then, by \eqref{firstregimeest},
$$kr\left|H^{(1)}_{\nu}(kr)\right|^2
\leq kr A_0^2\frac{1}{((kr)^2-\nu^2)^{1/2}}
\leq A_0^2\frac{1}{(1-(\nu/(kr))^2)^{1/2}}
\leq A_0^2\frac{\tilde{C}_0}{(\tilde{C}_0^2-1)^{1/2}}.
$$

Let us now assume that $j>\tilde{j}_0(kr)$, that is $kr< \tilde{C}_0\nu$.
First of all, by our more restrictive assumption on $B_0$ and by \eqref{extra},
we have that
$$kR_0\leq \tilde{C}_0\nu/B_0\leq 2\nu/(3\rme)
\leq
 \nu/2\leq \nu-C_0\nu^{1/3}.$$

 We have three different cases. In the first, we have that $kr\geq \nu+C_0\nu^{1/3}$.
 In the second,
 we have that $\nu+C_0\nu^{1/3}> kr> \nu-C_0\nu^{1/3}$.
 In the third case, $kr$ is less than or equal to $\nu-C_0\nu^{1/3}$.

In all cases we can use \eqref{thirdregimeest} for $z=kR_0$ and obtain that
$$A_0^{-1} \frac{1}{\sqrt{\nu}}\left(\frac{\nu}{\rme kR_0}\right)^{\nu}         \leq \left|H^{(1)}_{\nu}(kR_0)\right|.$$

In the first case, using \eqref{firstregimeest},
we have that
$$\left|H^{(1)}_{\nu}(kr)\right|\leq A_0\frac{1}{((kr)^2-\nu^2)^{1/4}}\leq
A_0\frac{1}{(2C_0\nu^{4/3}+C_0^2\nu^{2/3})^{1/4}}\leq \frac{A_0}{(2C_0)^{1/4}}\nu^{-1/3}.
$$
In the second case, by \eqref{secondregimeest}, we have that
$$\left|H^{(1)}_{\nu}(kr)\right|\leq A_0\nu^{-1/3},
$$
whereas in the third, using again \eqref{thirdregimeest}, we obtain
$$\left|H^{(1)}_{\nu}(kr)\right|\leq A_0
\left(\frac{2\nu}{kr}\right)^{\nu}
\frac{1}{(\nu^2-(kr)^2)^{1/4}}\exp\left(-(\nu^2-(kr)^2)^{1/2}\right).
$$
Since $(\nu^2-(kr)^2)\geq (\nu^2-(\nu-C_0\nu^{1/3})^2)\geq 
2C_0\nu^{4/3}-C_0^2\nu^{2/3}$, then by \eqref{extra2} we have that
$(\nu^2-(kr)^2)\geq C_0\nu^{4/3}\geq C_0z_0^{4/3}$. Hence
there exists a constant $\hat{C}_0\geq 1$, depending on $N$ only, such that
$$\left|H^{(1)}_{\nu}(kr)\right|\leq \hat{C}_0A_0
\left(\frac{2\nu}{kr}\right)^{\nu}.
$$

We conclude that in the first and second cases we have
\begin{equation}
\frac{\left|H^{(1)}_{\nu}(kr)\right|}{\left|H^{(1)}_{\nu}(kR_0)\right|}\leq \max\{1,1/(2C_0)^{1/4}\}
A_0^2\nu^{1/6}\left(\frac{\rme kR_0}{\nu}\right)^{\nu},
\end{equation}
whereas in the third we have
\begin{equation}
\frac{\left|H^{(1)}_{\nu}(kr)\right|}{\left|H^{(1)}_{\nu}(kR_0)\right|}\leq \hat{C}_0
A_0^2\sqrt{\nu}
\left(\frac{2\rme R_0}{r}\right)^{\nu}\leq
\hat{C}_0A_0^2\sqrt{\nu}
\left(\frac{2\rme }{B_0}\right)^{\nu}.
\end{equation}

Using the fact that $\nu\geq 2$ and our assumption on $B_0$, we conclude that for any
$j>\tilde{j}_0(kr)$ we have
\begin{equation}
\frac{\left|H^{(1)}_{\nu}(kr)\right|}{\left|H^{(1)}_{\nu}(kR_0)\right|}\leq \hat{A}_0
\sqrt{\nu}
\left(\frac{2}{3}\right)^{\nu}\leq 
\hat{A}_0\sqrt{\tilde{\nu}_0(kr)+1}
\left(\frac{2}{3}\right)^{\tilde{\nu}_0(kr)+1},
\end{equation}
where $\hat{A}_0=\max\{1,1/(2C_0)^{1/4}\}\hat{C}_0A_0^2$ is a constant depending on $N$ only.

Therefore, taking
$$\tilde{A}=\max\left\{(\pi/2)A_0^2\frac{\tilde{C}_0}{(\tilde{C}_0^2-1)^{1/2}},
\hat{A}_0^2
\right\},$$
and using \eqref{nf4} with
$j_0=\tilde{j}_0(kr)$ we have
\begin{equation}\label{nfll}
\|u^s\|^2_{L^2(\partial B_r)}\leq \tilde{A}B_1\left[\|u^s_{\infty}\|^2_{L^2(\mathbb{S}^{N-1})}
+
(\tilde{\nu}_0(kr)+1)
\left(\frac{2}{3}\right)^{2(\tilde{\nu}_0(kr)+1)}
\|u^s\|^2_{L^2(\partial B_{R_0})}
\right].
\end{equation}
Let us note that $\tilde{A}$ depends on $N$ only.
Moreover, $2(\tilde{\nu}_0(kr)+1)
\geq 2kr/\tilde{C}_0\geq 2z_0/\tilde{C}_0\geq 4$.
Hence
$$(\tilde{\nu}_0(kr)+1)
\left(\frac{2}{3}\right)^{2(\tilde{\nu}_0(kr)+1)}\leq 
(kr/\tilde{C}_0)
\left(\frac{2}{3}\right)^{2kr/\tilde{C}_0}.
$$
Finally, assuming that $k_1(\varepsilon,r)\geq k_0$, we obtain that
if $k\geq k_1(\varepsilon,r)$ then at least we have
\begin{multline}\label{nfagain}
\|u^s\|_{L^2(\partial B_r)}\leq\\ \sqrt{\tilde{A}B_1}\left[\|u^s_{\infty}\|^2_{L^2(\mathbb{S}^{N-1})}
+
(k_1(\varepsilon,r)r/\tilde{C}_0)\left(\frac{2}{3}\right)^{2k_1(\varepsilon,r)r/\tilde{C}_0}
\|u^s\|^2_{L^2(\partial B_{R_0})}
\right]^{1/2}.
\end{multline}

We summarize these results in the following proposition.

\begin{prop}\label{highfreqprop}

Let $N\geq 2$. We keep the previous notation and assumptions.
Let
$k_0>0$ be such that $k_0R_0\geq z_0$, $z_0$ as in Theorem~\textnormal{\ref{regimes}}
and such that \eqref{extra2} and \eqref{extra} are satisfied.

Let us assume that $B_0$ satisfies \eqref{B_0ass}.
Then there exists a constant $\tilde{A}$, depending on $N$ only, such that for any $k\geq k_0$ and any $r$, $B_0R_0\leq r\leq B_1R_0$, we have
\begin{equation}\label{nfllbis}
\|u^s\|_{L^2(\partial B_r)}\leq \sqrt{\tilde{A}B_1}\left[\varepsilon^2
+
M^2
(kr/\tilde{C}_0)
\left(\frac{2}{3}\right)^{2kr/\tilde{C}_0}
\right]^{1/2}.
\end{equation}

In particular, if $k\geq k_1(\varepsilon)$, where $k_1(\varepsilon)$ is as in Theorem~\textnormal{\ref{mainteo}} and satisfies $k_0\leq k_1(\varepsilon)$, we have that
\begin{equation}\label{nflltris}
\|u^s\|_{L^2(\partial B_r)}\leq \sqrt{\tilde{A}B_1}\left[\varepsilon^2
+
M^2
(k_1(\varepsilon)B_0R_0/\tilde{C}_0)
\left(\frac{2}{3}\right)^{2k_1(\varepsilon)B_0R_0/\tilde{C}_0}
\right]^{1/2}.
\end{equation}
\end{prop}

We notice that there exists an absolute constant $C_1\geq 1$
such that if $k_0\leq k_1(\varepsilon)$ we have
$$(k_1(\varepsilon)B_0R_0/\tilde{C}_0)
(2/3)^{2k_1(\varepsilon)B_0R_0/\tilde{C}_0}\leq C_1(3/4)^{2k_1(\varepsilon)B_0R_0/\tilde{C}_0}.$$
Hence \eqref{nflltris} may be written in the following simpler form
$$
\|u^s\|_{L^2(\partial B_r)}\leq \sqrt{\tilde{A}B_1}\left[\varepsilon^2
+
C_1M^{2(1-\tilde{\beta})}\varepsilon^{2\tilde{\beta}}
\right]^{1/2}\leq \sqrt{2C_1\tilde{A}B_1}
M^{(1-\tilde{\beta})}\varepsilon^{\tilde{\beta}}
$$
where we assume that $\varepsilon\leq M$ and
$$\tilde{\beta}=\frac{\log(4/3)}{\tilde{C}_0C_{\alpha}}\frac{B_0}{B_1}\left(\frac{2}{3\tilde{B}}\right)^{1+\alpha}\left(\frac{b_0\rme^{a_1}}{2}\right)^{\alpha}.$$

Finally, we wish to conclude with the following remark. Let us assume that for some positive exponent  $\tau$ and a constant $C_2$ we have, for any $k\geq k_0$,
\begin{equation}\label{Massum}
M\leq C_2k^{\tau}.
\end{equation}
Hence,  for any $r\geq B_0R_0$,
$$M\leq C_2k^{\tau}
\leq  C_2 (\tilde{C}_0/r)^{\tau}
(kr/\tilde{C}_0)^{\tau}\leq C_2 (\tilde{C}_0/B_0R_0)^{\tau} (kr/\tilde{C}_0)^{\tau}.
$$

Then there exists an absolute constant $C(\tau)$, depending on $\tau$ only, such that
\begin{multline*}
M^2(kr/\tilde{C}_0)
\left(\frac{2}{3}\right)^{2kr/\tilde{C}_0}\leq
C^2_2 (\tilde{C}_0/B_0R_0)^{2\tau}
C_3(\tau)\left(\frac{3}{4}\right)^{2kr/\tilde{C}_0}\leq\\
C^2_2 (\tilde{C}_0/B_0R_0)^{2\tau}
C_3(\tau)\left(\frac{3}{4}\right)^{2kB_0R_0/\tilde{C}_0}.
\end{multline*}
Therefore, the following corollary with a Lipschitz stability estimate holds.

\begin{cor}\label{Lipschitz}
Under the assumptions of Proposition~\textnormal{\ref{highfreqprop}}, let us further assume that, for some positive exponent  $\tau$ and a constant $C_2$, \eqref{Massum} holds for any $k\geq k_0$.

If $\varepsilon\leq 1/\rme$ and
$$k\geq \frac{\tilde{C}_0}{\log(4/3)B_0R_0}\log(1/\varepsilon)$$
we have
\begin{equation}\label{Lipest}
\|u^s\|_{L^2(\partial B_r)}\leq \left(
\tilde{A}B_1(1+C^2_2 (\tilde{C}_0/B_0R_0)^{2\tau}
C_3(\tau))\right)^{1/2}
\varepsilon.
\end{equation}
\end{cor}

\section{Stability estimates: from far-field up to the obstacle}\label{uptotheboundarysec}

We begin this section by establishing suitable a priori estimates for the solution to the direct scattering problem \eqref{Helmeq}.

We recall that, for any $k>0$, by $\Phi_{k}$ we denote the fundamental solution to the Helmholtz equation $\Delta u+k^2u =0$ which is given by
$$\Phi_k(x,y)=\frac{\rmi}{4}\left(\frac{k}{2\pi\|x-y\|}\right)^{(N-2)/2}H^{(1)}_{(N-2)/2}(k\|x-y\|)
\quad\text{for any }x,\ y\in\mathbb{R}^N,\ x\neq y.$$
We remark that for $N=2,3$ this reduces to the well known formulas
$$\Phi_k(x,y)=\frac{\rmi}{4}H^{(1)}_{0}(k\|x-y\|)\quad\text{for any }x,\ y\in\mathbb{R}^2,\ x\neq y.$$
and
$$\Phi_k(x,y)=\frac{\rme^{\rmi k\|x-y\|}}{4\pi\|x-y\|}\quad\text{for any }x,\ y\in\mathbb{R}^3,\ x\neq y.$$

Let now $\Sigma$ be a scatterer which may be characterized as the closure of a bounded Lipschitz open set. We recall that a bounded open set $D$ is said to be Lipschitz if 
for any $x\in\partial D$ there exist a positive $r$ and a Lipschitz function $\varphi:\mathbb{R}^{N-1}\to\mathbb{R}$, such that $\varphi(0)=0$ and, up to a rigid change of coordinates, we have $x=0$ and
$$B_r(x)\cap D = \{y=(y',y_N)\in B_r(x):\ y_N<\varphi(y')\}.$$

For any density $\psi\in L^2(\partial\Sigma)$, let us define $w$ and $v$ as the corresponding
\emph{single-layer} and \emph{double-layer potentials} with density $\psi$, namely
\begin{equation}\label{single}
w(x)=\int_{\partial\Sigma}\psi(y)\Phi_k(x,y)\rmd\mathcal{H}^{N-1}(y)\quad x\in \mathbb{R}^N\backslash \Sigma,
\end{equation}
and
\begin{equation}\label{double}
v(x)=\int_{\partial\Sigma}\psi(y)\frac{\partial \Phi_k(x,y)}{\partial\nu(y)}\rmd\mathcal{H}^{N-1}(y)\quad x\in \mathbb{R}^N\backslash \Sigma,
\end{equation}
where $\nu$ is the exterior normal to $\Sigma$.
We notice that $w$ and $v$ satisfy the Helmholtz equation as well as the Sommerfeld radiation condition. We call $w_{\infty}$ and $v_{\infty}$ their far-field patterns, respectively.
Then the following result holds.


\begin{prop}\label{apririestprop}
Under the previous notation and assumptions we have the following estimates.
For any $x\in G=\mathbb{R}^N\backslash\Sigma$ let $d=\mathrm{dist}(x,\Sigma)$. Then there exists a constant $C$, depending on $N$ only, such that for any $k>0$ and any $x\in G$ we have
\begin{equation}\label{estimate13}
|w(x)|\leq C(\mathcal{H}^{N-1}(\partial\Sigma))^{1/2}\|\psi\|_{L^2(\partial\Sigma)}
\frac{1}{d^{N-2}}\max\{1,(kd)^{(N-3)/2}\}\quad \text{for }N\geq 3
\end{equation}
and
\begin{equation}\label{estimate12}
|w(x)|\leq (1/4)(\mathcal{H}^{N-1}(\partial\Sigma))^{1/2}\|\psi\|_{L^2(\partial\Sigma)}
|H^{(1)}_0(kd)|\quad \text{for }N=2,
\end{equation}
and
\begin{equation}\label{estimate2}
|v(x)|\leq C(\mathcal{H}^{N-1}(\partial\Sigma))^{1/2}\|\psi\|_{L^2(\partial\Sigma)}
\frac{1}{d^{N-1}}\max\{1,(kd)^{(N-1)/2}\}\quad \text{for }N\geq 2.
\end{equation}

Furthermore, for any $N\geq 2$, any $k>0$ and any $\hat{x}\in \mathbb{S}^{N-1}$, we have
\begin{equation}\label{sinfty}
w_{\infty}(\hat{x})=\frac{\rmi}{2}\frac{\rme^{-(N-1)\pi\rmi/4}}{(2\pi)^{(N-1)/2}}k^{(N-3)/2}
\int_{\partial\Sigma}\psi(y)\rme^{-\rmi k \hat{x}\cdot y}
\mathrm{d}\mathcal{H}^{N-1}(y)
\end{equation}
and
\begin{equation}\label{dinfty}
v_{\infty}(\hat{x})=\frac{\rmi}{2}\frac{\rme^{-(N-1)\pi\rmi/4}}{(2\pi)^{(N-1)/2}}k^{(N-3)/2}\int_{\partial\Sigma}\psi(y)\frac{\partial\rme^{-\rmi k \hat{x}\cdot y}}{\partial \nu(y)}
\mathrm{d}\mathcal{H}^{N-1}(y).
\end{equation}
Therefore, for any $N\geq 2$, any $k>0$ and any $\hat{x}\in \mathbb{S}^{N-1}$, we have
\begin{equation}\label{sfinal}
|w_{\infty}(\hat{x})|\leq \frac{k^{(N-3)/2}}{2(2\pi)^{(N-1)/2}}(\mathcal{H}^{N-1}(\partial\Sigma))^{1/2}\|\psi\|_{L^2(\partial\Sigma)}
\end{equation}
and
\begin{equation}\label{dfinal}
|v_{\infty}(\hat{x})|\leq \frac{k^{(N-1)/2}}{2(2\pi)^{(N-1)/2}}(\mathcal{H}^{N-1}(\partial\Sigma))^{1/2}\|\psi\|_{L^2(\partial\Sigma)}.
\end{equation}
\end{prop}

\proof{.}
Let us notice that for any $x$, $y\in\mathbb{R}^N$, $x\neq y$, we have 
\begin{multline*}
\nabla_y\Phi_k(x,y)=\\
\frac{\rmi}{4}\frac{k^{N-2}k}{(2\pi)^{(N-2)/2}}\left[\frac{H^{(1)}_{(N-4)/2}(k\|x-y\|)}{(k\|x-y\|)^{(N-2)/2}}-(N-2)\frac{H^{(1)}_{(N-2)/2}(k\|x-y\|)}{(k\|x-y\|)^{N/2}}\right]\frac{y-x}{\|y-x\|}.
\end{multline*}
We also remark that for $N=2,3$ this reduces to
$$\nabla_y\Phi_k(x,y)=
-\frac{\rmi}{4}k
\left[H^{(1)}_{1}(k\|x-y\|)\right]\frac{y-x}{\|y-x\|}\quad N=2$$
and
$$\nabla_y\Phi_k(x,y)=
\frac{\rmi}{4\sqrt{2\pi}}k^2
\left[\frac{\rmi H^{(1)}_{1/2}(k\|x-y\|)}{(k\|x-y\|)^{1/2}}-\frac{H^{(1)}_{1/2}(k\|x-y\|)}{(k\|x-y\|)^{3/2}}\right]\frac{y-x}{\|y-x\|}\quad N=3.$$

Then the estimates \eqref{estimate13} and \eqref{estimate2} follow by straightforward, although lengthy, computations. The main ingredient is the asymptotic behavior of Hankel functions
$H_{\nu}(z)$, with $\nu\geq 0$, as $z\to 0^+$ and $z\to +\infty$. The latter is given in \eqref{argument}, while the former is the following
\begin{equation}\label{rto0}
H^{(1)}_{\nu}(z)\sim \left\{\begin{array}{ll}
\vphantom{\Bigg(}\displaystyle{-\rmi\frac{2}{\pi}\log(2/z)}
& \text{for }\nu=0\\
\vphantom{\Bigg(}\displaystyle{-\rmi\frac{\Gamma(\nu)}{\pi}(2/z)^{\nu}}& \text{for }\nu>0
\end{array}\right.
\quad\text{as }z\to 0^+.
\end{equation}

For what concerns \eqref{estimate12}, the estimate follows trivially from this remark. We have that 
$|H^{(1)}_{\nu}(z)|^2$ is a decreasing function of $z>0$ for any fixed $\nu\geq 0$, see for instance \cite[page~446]{Wat}. This remark and the fact that
$$|H^{(1)}_0(z)|\sim \frac{2}{\pi}\log(2/z)\text{ as }z\to 0^+\quad\text{and}\quad |H^{(1)}_0(z)|\sim\left(\frac{2}{\pi z}\right)^{1/2}\text{ as }z\to +\infty$$
provide also the correct way to interpret \eqref{estimate12}.

The relationships \eqref{sinfty} and \eqref{dinfty} may be proved, for $N=3$, by the argument used to prove Theorem~2.5 in \cite{Col e Kre98}. For $N\neq 3$, a standard modification is needed.

Finally, \eqref{sfinal} and \eqref{dfinal} follow immediately from \eqref{sinfty} and \eqref{dinfty}, respectively.\cvd

\bigskip

We also need the following easy lemma.

\begin{lem}\label{cacciolemma}
Let $\Sigma\subset \overline{B}_R$, $R>0$, be any scatterer, with no regularity assumption.
Let us fix $k>0$ and $\omega\in\mathbb{S}^{N-1}$. Let $R<r_1<r$.
We let $u=u(\omega,k,\Sigma)$ be the solution to \eqref{Helmeq}.
Then there exists a constant $C$, depending on $R$, $r_1$ and $r$ only, such that
$$\|u\|_{H^1(B_{r_1}\backslash\Sigma)}\leq C\max\{1,k\}\|u\|_{L^2(B_{r}\backslash\Sigma)}.$$
\end{lem}

\proof{.}
We fix $\tilde{r}$ such that $R<r_1<\tilde{r}<r$, $\tilde{r}$ depending on $R$, $r_1$ and $r$ only. 
Then
by a standard Caccioppoli's inequality, we infer that
$$\|u\|_{H^1(B_{\tilde{r}}\backslash\overline{B_{r_1}})}\leq C_1\max\{1,k\}\|u\|_{L^2(B_{r}\backslash\Sigma)},$$
for some constant $C_1$ depending on $R$, $r_1$ and $r$ only.

By integrating over $B_{\tilde{r}}\backslash\overline{B_{r_1}}$ in spherical coordinates,
we infer that there exists $\rho$, $r_1<\rho<\tilde{r}$, such that
$$\int_{\partial B_{\rho}}|u|^2\leq \frac{3}{\tilde{r}-r}_1\|u\|^2_{L^2(B_{r}\backslash\Sigma)}
\quad\text{and}\quad\int_{\partial B_{\rho}}\left|\frac{\partial u}{\partial\nu}\right|^2\leq \frac{3}{\tilde{r}-r_1}C_1^2\max\{1,k^2\}\|u\|^2_{L^2(B_{r}\backslash\Sigma)}.$$

Since
$$\int_{B_{\rho}\backslash\Sigma}\|\nabla u\|^2=k^2\int_{B_{\rho}\backslash\Sigma}|u|^2+\int_{\partial B_{\rho}}\frac{\partial u}{\partial\nu}\overline{u},$$
the thesis immediately follows.\cvd

\bigskip

Let us remark that, for the same ${\rho}$ as in the previous proof,
if $u^s=u^s(\omega,k,\Sigma)$ we have

$$\|u^s\|_{L^2(\partial B_{\rho})}\leq \left(\frac{3}{\tilde{r}-r_1}\right)^{1/2}\|u\|_{L^2(B_{r}\backslash\Sigma)}+(\mathcal{H}^{N-1}(\partial B_{\rho}))^{1/2}$$
and
$$\left\|\frac{\partial u^s}{\partial\nu}\right\|_{L^2(\partial B_{\rho})}\leq \left(\frac{3}{\tilde{r}-r_1}\right)^{1/2}C_1\max\{1,k\}\|u\|_{L^2(B_{r}\backslash\Sigma)}+k(\mathcal{H}^{N-1}(\partial B_{\rho}))^{1/2}.$$
Therefore, since $u^s$ has the following Helmholtz representation
\begin{equation}\label{Helmholtzrepr}
u^s(x)=\int_{\partial B_{\rho}}
\frac{\partial u^s(y)}{\partial\nu}\Phi_k(x,y)-u^s(y)\frac{\partial \Phi_k(x,y)}{\partial\nu(y)}
\rmd\mathcal{H}^{N-1}(y)
\quad \|x\|>\rho,
\end{equation}
we can use all the results given in Proposition~\ref{apririestprop}.
In particular we have for any $k>0$
\begin{multline}\label{apriorifarfield}
\|\mathcal{A}(\Sigma)(\cdot,\omega,k)\|_{L^2(\mathbb{S}^{N-1})}\leq\\
\mathcal{H}^{N-1}(\mathbb{S}^{N-1})
\frac{k^{(N-3)/2}}{2(2\pi)^{(N-1)/2}}\tilde{r}^{(N-1)/2}
\left[k
\|u^s\|_{L^2(\partial B_{\rho})}
+\left\|\frac{\partial u^s}{\partial\nu}\right\|_{L^2(\partial B_{\rho})}\right]\leq\\
C(N)
\left(\frac{k\tilde{r}}{2\pi}\right)^{(N-1)/2}
\left[\left(\frac{3}{\tilde{r}-r_1}\right)^{1/2}(1+C_1\max\{k^{-1},1\})\|u\|_{L^2(B_{r}\backslash\Sigma)}+2(\mathcal{H}^{N-1}(\partial B_{\tilde{r}}))^{1/2}\right]
\end{multline}
where $C(N)=\mathcal{H}^{N-1}(\mathbb{S}^{N-1})/2$.
Furthermore, by \eqref{Nedelec},
we also have for any $k>0$
\begin{multline}\label{Ned2}
\|\mathcal{A}(\Sigma)(\cdot,\omega,k)\|^2_{L^2(\mathbb{S}^{N-1})}=2\left(\frac{2\pi}{k}\right)^{(N-1)/2}
\Im\left(
\rme^{(N-3)\pi\rmi/4}
\mathcal{A}(\Sigma)(\omega,\omega,k)
\right)\leq\\
k^{-1}
(\mathcal{H}^{N-1}(\partial B_{\rho}))^{1/2}\left[k
\|u^s\|_{L^2(\partial B_{\rho})}
+\left\|\frac{\partial u^s}{\partial\nu}\right\|_{L^2(\partial B_{\rho})}\right]\leq\\
(\mathcal{H}^{N-1}(\partial B_{\tilde{r}}))^{1/2}\left[\left(\frac{3}{\tilde{r}-r_1}\right)^{1/2}(1+C_1\max\{k^{-1},1)\}\|u\|_{L^2(B_{r}\backslash\Sigma)}+2(\mathcal{H}^{N-1}(\partial B_{\tilde{r}}))^{1/2}\right].
\end{multline}
We can conclude that there exists a constant $C$, depending on $N$,
$R$, $r$ and $\|u\|_{L^2(B_{r}\backslash\Sigma)}$ only,
such that
\begin{equation}\label{aprioriffN3}
\|\mathcal{A}(\Sigma)(\cdot,\omega,k)\|_{L^2(\mathbb{S}^{N-1})}\leq C
\quad \text{for any }N\geq 3\text{ and any }k>0
\end{equation}
and
\begin{equation}\label{aprioriffN2}
\|\mathcal{A}(\Sigma)(\cdot,\omega,k)\|_{L^2(\mathbb{S}^{N-1})}\leq C\max\{k^{-1/2},1\}
\quad \text{for }N= 2\text{ and any }k>0.
\end{equation}

We continue by establishing suitable a priori estimates on the solutions to direct scattering problems with sound-soft scatterers. In particular we are interested in the high frequency asymptotics. The main tool will be provided by the results established by Chandler-Wilde and Monk, \cite{Cha-Mon}.

Let us fix an integer $m\geq 1$ and positive constants $\beta$, $R_0$, and $\delta_0$
such that
$R_0<R_0+\delta_0=R\leq \beta$. Let us also fix $r>R$.

We begin with the case in which $k$ lies in a bounded interval, that is we fix constants $0<k_0<k_1$ and we denote, for any $N\geq 2$,
\begin{equation}\label{I_Ndefin}
I_N=\left\{\begin{array}{ll}
[k_0,k_1]&\text{if }N=2,\\
(0,k_1]&\text{if }N\geq 3.
\end{array}\right.
\end{equation}

We have the following a priori estimate which follows from arguments used in \cite{Isak92}. We notice that a much more general class of scatterers may be used, see for instance \cite{Ron03} for sound-soft scatterers and \cite{Men-Ron} for the corresponding sound-hard case.

\begin{prop}\label{kboundedprop}
Under the previous notation, there exists a constant $C$, depending on $N$, $m$, $\beta$, $R_0$, $R$, $r$ and $I_N$ only, such that for any $k\in I_N$, any $\omega\in \mathbb{S}^{N-1}$ and any $\Sigma\in X(m,\beta,R_0,\delta_0)$ we have
\begin{equation}\label{aprioriL2}
\|u(\omega,k,\Sigma)\|_{L^2(B_r\backslash\Sigma)}\leq C.
\end{equation}
\end{prop}

We now consider the high frequency asymptotics. We limit ourselves to the case $m\geq 2$ and $N=2,3$. Fix $\Sigma\in X(m,\beta,R_0,\delta_0)$, $k>0$ and $\omega\in\mathbb{S}^{N-1}$.
We have that $u^s=u^s(\omega,k,\Sigma)$, the scattered field of the solution to the direct scattering problem \eqref{Helmeq}, may be described as the sum of a double- and a  single-layer potential in the following way
\begin{equation}\label{Helrep}
u^s(x)=\int_{\partial\Sigma}\psi(y)\left[\frac{\partial \Phi_k(x,y)}{\partial\nu(y)}-\rmi k\Phi_k(x,y)\right]\rmd\mathcal{H}^{N-1}(y)\quad x\in \mathbb{R}^N\backslash \Sigma,
\end{equation}
where $\psi\in C^0(\partial\Sigma)$, see for instance \cite{Col e Kre98}. 
Here $\psi$ solves the following integral equation
$$(I+K_k-\rmi k S_k)\psi=A\psi=-2 u^i\quad\text{on }\partial\Sigma$$
where $K_k$ and $S_k$ are the double- and single-layer operators, respectively, 
defined by, see Chapter~3 of \cite{Col e Kre98},
$$ K_k(\psi)(x)=\int_{\partial \Sigma} \frac{\partial \Phi_k(x,y)}{\partial\nu(y)} \psi(y)\rmd \mathcal{H}^{N-1}(y),\quad S_k(\psi)(x)=\int_{\partial \Sigma} \Phi_k(x,y) \psi(y)\rmd \mathcal{H}^{N-1}(y),$$
for almost every $x\in \partial \Sigma$. We also denote 
for any $x\in \mathbb{R}^N\backslash \partial \Sigma$
$$ \tilde{K}_k(\psi)(x)=\int_{\partial \Sigma} \frac{\partial \Phi_k(x,y)}{\partial\nu(y)} \psi(y)\rmd \mathcal{H}^{N-1}(y),\quad \tilde{S}_k(\psi)(x)=\int_{\partial \Sigma} \Phi_k(x,y) \psi(y)\rmd \mathcal{H}^{N-1}(y).$$

We have that $A=A(k,\Sigma)$ is bounded and bijective from $C^0(\partial\Sigma)$ onto itself and also from $L^2(\partial\Sigma)$ onto itself.

In the following theorem, we state some useful estimates on $A(k,\Sigma)^{-1}$,
$\tilde{K}_k$ and $\tilde{S}_k$.


\begin{teo}\label{Cha-Monteo}
Under the previous notation,
let $m\geq 2$ and $N=2,3$.
There exist constants $\tilde{C}$, depending on $N$, $m$, $\beta$ and $R_0$ only, and $\tilde{C}_1$, depending on $N$, $m$, $\beta$, $R_0$ and $r$ only,
such that 
for any $k$ satisfying $kR_0\geq 1$ and any $\Sigma\in X(m,\beta,R_0,\delta_0)$ we have

\begin{equation}\label{A-1}
\|A(k,\Sigma)^{-1}\|_{\mathcal{L}(L^2(\partial\Sigma),L^2(\partial\Sigma))}\leq \tilde{C}.
\end{equation}

Furthermore we also have
\begin{equation}\label{K-k}
\|\tilde{K}_k\|_{\mathcal{L}(L^2(\partial\Sigma),L^2(B_r \backslash \Sigma))}\leq \tilde{C}_1 k
\end{equation}
and
\begin{equation}\label{S-k}
\|\tilde{S}_k\|_{\mathcal{L}(L^2(\partial\Sigma),L^2(B_r \backslash \Sigma))}\leq \tilde{C}_1.
\end{equation}
\end{teo}

\proof{.}
The estimate \eqref{A-1} is an immediate corollary of Theorem~4.3 and Corollary~4.4 in \cite{Cha-Mon}, whereas 
\eqref{K-k} follows
from Theorem~5.2 in \cite{Mel}.

For what concerns the estimate on $\tilde{S}_k$, the argument depends on the dimension $N$. For $N=3$ we have that
$|\Phi_k(x,y)|\leq \Phi_0(x,y)=1/(4\pi\|x-y\|)$, therefore
$$\|\tilde{S}_k\|_{\mathcal{L}(L^2(\partial\Sigma),L^2(B_r \backslash \Sigma))} \leq \|\tilde{S}_0\|_{\mathcal{L}(L^2(\partial\Sigma),L^2(B_r \backslash \Sigma))}$$
where $\tilde{S}_0$ is the corresponding operator with $\Phi_k$ replaced by $\Phi_0$,
the fundamental solution of the Laplacian for $N=3$.
For $N=2$ the argument is slightly more involved. We recall that
$|H^{(1)}_0(z)|$ is a decreasing function of $z>0$, therefore, if we set for the time being $k_0=1/R_0$, for any $k\geq k_0$ we have $|\Phi_k(x,y)|\leq |\Phi_{k_0}(x,y)|$. Furthermore, there exist positive constants $C_1$ and $C_2$, depending on $R_0$ and $r$ only, such that
for any $x,y\in B_r\subset\mathbb{R}^2$ we have
 $$|\Phi_{k_0}(x,y)|\leq C_1\Phi_0(x,y)+C_2$$
where $\Phi_0(x,y)=-(2\pi)^{-1}\log(\|x-y\|)$ is the fundamental solution of the Laplacian for $N=2$. Here we have made use of \eqref{rto0}.
Therefore, for any $x\in B_r\backslash \Sigma$ and any $k\geq k_0$, we have 
$$|\tilde{S}_k(\psi)(x)|\leq \int_{\partial \Sigma} |\Phi_{k_0}(x,y)| |\psi(y)|\rmd \mathcal{H}^{N-1}(y)\leq C_1\tilde{S}_0(|\psi|)(x)+C_2E^{1/2}\|\psi\|_{L^2(\partial\Sigma)},$$
 $E$ as in \eqref{areabound}, and the conclusion immediately follows.\cvd

\bigskip

We may conclude with our a priori estimates stating the following result.


\begin{cor}\label{aprioricorollary}
Under the previous notation, let $m\geq 2$ and $N=2,3$.

There exists a constant $C_1$, depending on $N$, $m$, $\beta$, $R_0$, $\delta_0$ and $r$ only, such that 
for any $k$ satisfying $kR_0\geq 1$, any $\omega\in \mathbb{S}^{N-1}$ and any $\Sigma\in X(m,\beta,R_0,\delta_0)$ we have the following estimates.
We let $u=u(\omega,k,\Sigma)$, $u^s=u^s(\omega,k,\Sigma)$ and
$u^s_{\infty}=u^s_{\infty}(\omega,k,\Sigma)$. Then

\begin{equation}\label{l2estimate}
\|u(\omega,k,\Sigma)\|_{L^2(B_r\backslash\Sigma)}\leq C_1k
\end{equation}
and, for any $x\in G=\mathbb{R}^N\backslash\Sigma$, if $d=\mathrm{dist}(x,\Sigma)$ we have
\begin{equation}
|u^s(x)|\leq 4C\tilde{C}E
\left[
\max\left\{\frac{1}{d^{N-1}},\left(\frac{k}{d}\right)^{(N-1)/2}\right\}
\right]\quad\text{for }N=3,
\end{equation}
and 
\begin{equation}
|u^s(x)|\leq 2\tilde{C}E
\left[(1/4)k|H^{(1)}_0(kd)|+C\max\left\{\frac{1}{d^{N-1}},\left(\frac{k}{d}\right)^{(N-1)/2}\right\}
\right]\quad\text{for }N=2,
\end{equation}
where $C$ is as in Proposition~\textnormal{\ref{apririestprop}}. Notice that, if we fix $d_0>0$, then
$$k|H^{(1)}_0(kd)|\leq \hat{C}\sqrt{\frac{k}{d}}\leq \hat{C}\max\left\{\frac{1}{d},\sqrt{\frac{k}{d}}\right\}\quad \text{for any }d\geq d_0$$
where $\hat{C}$ depends on $d_0/R_0$ only.
Finally, we have for any $\hat{x}\in\mathbb{S}^{N-1}$
\begin{equation}
|u^s_{\infty}(\hat{x})|\leq 2\frac{k^{(N-1)/2}}{(2\pi)^{(N-1)/2}}\tilde{C}E
\end{equation}
hence by \eqref{Nedelec}
\begin{equation}\label{ff3}
\|u^s_{\infty}\|^2_{L^2(\mathbb{S}^{N-1})}\leq 
4\tilde{C}E.
\end{equation}
\end{cor}

\proof{.} We sketch the proof only of estimate \eqref{l2estimate}, all the others follow in a standard way from the previously stated results.

We notice that if
$\psi=A(k,\Sigma)^{-1}(-2u^i)$, then
\begin{equation}\label{psi}
\|\psi\|_{L^2(\partial\Sigma)}\leq 2\tilde{C}(\mathcal{H}^{N-1}(\partial\Sigma))^{1/2}\leq 2\tilde{C}E^{1/2},
\end{equation}
where $\tilde{C}$ is as in \eqref{A-1} and $E$ is as in \eqref{areabound}.

From \eqref{Helrep} we have
$$
\|u^s(\omega,k,\Sigma)\|_{L^2(B_r\backslash\Sigma)} \leq \left[\|\tilde{K}_k\|_{\mathcal{L}(L^2(\partial\Sigma),L^2(B_r \backslash \Sigma))} +k \|\tilde{S}_k\|_{\mathcal{L}(L^2(\partial\Sigma),L^2(B_r \backslash \Sigma))} \right]
\|\psi\|_{L^2(\partial\Sigma)}.
$$
Recalling that $\|u^i\|_{L^2(B_r\backslash\Sigma)} \leq |B_1|^{1/2}r^{N/2}$, the estimate \eqref{l2estimate} follows then from \eqref{K-k}, \eqref{S-k} and
\eqref{psi}.\cvd

\bigskip

We remark that with the use of Corollary~\ref{aprioricorollary} and of the results obtained
by coupling
Proposition~\ref{apririestprop},
Lemma~\ref{cacciolemma} and the following discussion, in particular \eqref{aprioriffN3} and \eqref{aprioriffN2}, with Proposition~\ref{kboundedprop}, we may obtain a priori estimates for solutions of the direct scattering problem \eqref{Helmeq}  for any $k>0$ if $N=3$ and for any $k\geq k_0$ for $N=2$, $k_0$ being a fixed positive constant. We also notice that an estimate related to \eqref{ff3} in the high frequencies regime was obtained in \cite{Lak-Vai}. 


Under the previous notation, let $m\geq 1$ and $N\geq 2$.
We consider $\Sigma_1$, $\Sigma_2\in X(m,\beta,R_0,\delta_0)$ and we call $K$ the convex hull of $\Sigma_1\cup\Sigma_2$.
We begin by observing that $K$ is closed, convex with not empty interior.
Moreover, $B_{R_0}\subset K\subset \overline{B_R}$.

We fix $k>0$ and $\omega\in\mathbb{S}^{N-1}$. Let $u_1=u(\omega,k,\Sigma_1)$
and $u_2=u(\omega,k,\Sigma_2)$
be the solutions to \eqref{Helmeq} with $\Sigma$ replaced by $\Sigma_1$ and $\Sigma_2$ respectively. Let the corresponding scattered fields be denoted by
$u^s_1=u^s(\omega,k,\Sigma_1)$
and $u^s_2=u^s(\omega,k,\Sigma_2)$ and their far-field patterns by
$\mathcal{A}(\Sigma_1)(\cdot,\omega,k)$ and $\mathcal{A}(\Sigma_2)(\cdot,\omega,k)$
 respectively. We call $v=u_1-u_2=u_1^s-u_2^s$. Obviously, $v=v(\omega,k)$.

We wish to estimate $v$ outside $K$ in terms of the difference between the far-field patterns
$\mathcal{A}(\Sigma_1)(\cdot,\omega,k)$ and $\mathcal{A}(\Sigma_2)(\cdot,\omega,k)$.
We are clearly interested only in the high frequencies regime. Therefore,
let $z_0$, $C_0$ and $A_0$ be as in Theorem~\ref{regimes} and such that
\eqref{extra2} and \eqref{extra} are satisfied.
We fix a constant $k_0\geq 1$ such that $k_0R_0\geq z_0$. We fix $\rho_1=\tilde{B}_0(R+1)$,
where $\tilde{B}_0$ is as in \eqref{B_0ass}.

We begin with the following remark.
There exists $\rho_0$, $R<\rho_0<R+1$, such that
\begin{equation}
\int_{\partial B_{\rho_0}}|v|^2\leq \int_{B_{R+1}\backslash \overline{B_{R}}}|v|^2.
\end{equation}
Then $\rho_1=B_0\rho_0$ where $\tilde{B}_0\leq B_0\leq \tilde{B}_0(R+1)/R$.

We assume that $k\geq k_0$ and that for some $\varepsilon>0$ and $M>0$ we have
\begin{equation}\label{apriori}
\|\mathcal{A}(\Sigma_1)(\cdot,\omega,k)-\mathcal{A}(\Sigma_2)(\cdot,\omega,k)\|_{L^2(\mathbb{S}^{N-1})}\leq\varepsilon,\quad\|v\|_{L^2(B_{R+1})\backslash \overline{B_R}}\leq M.
\end{equation}
and we assume that $\varepsilon\leq M$.
Then for any $\tilde{B}_1\geq\tilde{B}_0$ we set
$\rho_2=\tilde{B}_1(R+1)=B_1\rho_0$ where $\tilde{B}_1\leq B_1\leq \tilde{B}_1(R+1)/R$.
We assume, for the time being, that
\begin{equation}\label{firststepest}
\|v\|_{L^2(B_{\rho_2}\backslash\overline{B_{\rho_1}})}\leq \eta_1(\varepsilon,k,M,\tilde{B}_1).
\end{equation}
We note that
the results of the previous section allow us to estimate precisely
$\eta_1(\varepsilon,k,M,\tilde{B}_1)$. Our aim is to estimate $v$ up to the boundary of $K$ and precisely on the following set $B_{\rho_1}\backslash K$. We shall use the results in
\cite{Sub-Isa}, a consequence of those in \cite{Hry-Isa}.

We begin with the following intermediate case. Let us take $P\in \partial K$. Let $\pi$ be a supporting hyperplane for $K$ passing through $P$ and let $S^+$ be the open half-space with boundary $\pi$ not intersecting $K$. We wish to estimate the $L^2$ norm of $v$ on
$B_{\rho_1}\cap S^+$. We argue in the following way. Let $\nu$ be the normal to $\pi$ pointing inside $S^+$ and let $P_1$ be the point of $\pi$ intersecting the half-line $l=\{s\nu:\ s\geq 0\}$.
We notice that $0< R_0\leq \|P_1\|\leq R$.
For any fixed $r$, $\rho_1+1\leq r\leq \rho_1+2$, we consider the open cylinder $T_r$, contained in $S^+$, whose lower base is contained in
$\pi$, is centered in $P_1$ and has radius $r$, and whose height is $r$. We call $\Gamma_r$ 
its upper base, that is the one contained in $S^+$. Let us now fix 
\begin{equation}
\tilde{B}_1=2\tilde{B}_0+5\text{ and }\rho_2=\tilde{B}_1(R+1).
\end{equation}
We notice that for any
$r$, $\rho_1+1\leq r\leq \rho_1+2$, we have
$$B_{\rho_1}\cap S^+\subset T_r\subset B_{\rho_2-1}\cap S^+
\quad\text{and}\quad
\Gamma_r\subset B_{\rho_2-1}\backslash \overline{B_{\rho_1+1}}.$$

Then, by the same argument used in Lemma~\ref{cacciolemma},
we infer that there exists $r$, $\rho_1+1< r< \rho_1+2$, such that
\begin{equation}\label{Cauchy_data}
\int_{\Gamma_r}|v|^2\leq 3
\|v\|^2_{L^2(B_{\rho_2}\backslash \overline{B_{\rho_1}})}\quad\text{and}\quad
\int_{\Gamma_r}|\nabla v|^2\leq 3C_1^2k^2\|v\|^2_{L^2(B_{\rho_2}\backslash \overline{B_{\rho_1}})}
\end{equation}
for some constant $C_1$ depending on $N$, $R_0$ and $\delta_0$ only. We conclude that
\begin{equation}\label{Cauchy_data2}
\|v\|_{L^2(\Gamma_r)}+\|\nabla v\|_{L^2(\Gamma_r)}\leq
\sqrt{3}(1+C_1k)
\|v\|_{L^2}(B_{\rho_2}\backslash \overline{B_{\rho_1}})
\leq C_2k
\eta_1(\varepsilon,k,M,\tilde{B}_1)
\end{equation}
where $C_2=\sqrt{3}(1+C_1)$.

Finally, let us assume that for some constant $\tilde{M}$, $\varepsilon\leq \tilde{M}$,
we have
\begin{equation}\label{aprioridataL2}
\|u_1\|_{L^2(B_{\rho_2}\backslash \Sigma_1)},\ \|u_2\|_{L^2(B_{\rho_2}\backslash \Sigma_2)}\leq \tilde{M}/2
\end{equation}
hence, again by Lemma~\ref{cacciolemma}, there exists a constant $C_3$, depending on $N$, $R_0$ and $\delta_0$ only, such that
\begin{equation}\label{aprioridataH1}
\|v\|_{H^1(T_r)}\leq C_3k\tilde{M}.
\end{equation}
Let us also notice that we can assume
$$M\leq \tilde{M}.$$

We obtain the following lemma, an immediate consequence of Theorem~1.1 in \cite{Sub-Isa}.

\begin{lem}\label{hyperplanelemma}
Under the previous notation and assumptions,
let us assume that $\varepsilon\leq M\leq \tilde{M}$ and that
$$
\eta_1=\eta_1(\varepsilon,k,\tilde{M},\tilde{B}_1)\leq (C_3/C_2)\tilde{M}.
$$

Then we have that for any $k\geq k_0$
\begin{equation}\label{L2halfestimate}
\|v\|_{L^2(B_{\rho_1}\cap S^+)}\leq C_4\eta_2(\varepsilon,k,\tilde{M})
\end{equation}
where $C_4$ depends on $N$, $R_0$ and $\delta_0$ only and
\begin{equation}\label{eta_2defin}
\eta_2(\varepsilon,k,\tilde{M})=
\left(
C_2^2k^2
\eta^2_1+\frac{
C_3^2k^2
\tilde{M}^2}{\left(-\log\left((C_2/C_3)\eta_1/\tilde{M}\right)+k\right)^{1/8}}
\right)^{1/2}.
\end{equation}
\end{lem}

To proceed further we need to make the following additional assumption, namely that $m\geq 2$.
Let us begin by studying some geometrical properties of $K$,  the convex hull of $\Sigma_1\cup\Sigma_2$. Let $P$ be any point belonging to $\partial K$.
Without loss of generality, let us assume that $P=s_0e_N$, where $s_0>0$ and $e_1,\ldots,e_N$ denote the canonical base in $\mathbb{R}^N$. Then there exist constants $r_0>0$ and $\theta_0$,
$0<\theta_0<\pi/2$, depending on $m$, $\beta$, $R_0$ and $\delta_0$ only, such that the following holds. There exists a point $Q$, depending on $P$, such that $\|P-Q\|=r_0$,
$B_{r_0}(Q)\subset K$ and finally the angle between $-e_N$ and the vector $Q-P$ is at most $\theta_0$. 
We immediately infer a few interesting properties. 
First of all, there exists a
unique supporting hyperplane for $K$ passing through $P$, the hyperplane whose normal is given by the vector $Q-P$.
For any direction $\hat{x}\in\mathbb{S}^{N-1}$ there exists a unique $s_0(\hat{x})>0$ such that
$s\hat{x}\in K$ for any $0\leq s\leq s_0(\hat{x})$ and $s\hat{x}\not\in K$ for any $s>s_0(\hat{x})$. Clearly $P(\hat{x})=s_0(\hat{x})\hat{x}$ is the only one of these points belonging to $\partial K$ and $R_0\leq s_0(\hat{x})\leq R$ for any $\hat{x}\in\mathbb{S}^{N-1}$.
For any $\hat{x}\in\mathbb{S}^{N-1}$ we denote by $\pi(\hat{x})$ the unique supporting hyperplane for $K$ passing through $P(\hat{x})$ and by $S^+(\hat{x})$ the open half-space with 
boundary $\pi(\hat{x})$ not intersecting $K$.
A further crucial geometrical property of $K$ is given in the following lemma.

\begin{lem}\label{Kproperties}
Let $m\geq 2$ and $N\geq 2$. Let us consider $\Sigma_1$, $\Sigma_2\in X(m,\beta,R_0,\delta_0)$ and let $K$ be the convex hull of $\Sigma_1\cup\Sigma_2$.

Let $\hat{x}\in\mathbb{S}^{N-1}$ and let $P=P(\hat{x})=s_0(\hat{x})\hat{x}$
belong to $\partial K$.
For any $d>0$ we denote $P_d=(s_0(\hat{x})+d)\hat{x}$.

Then there exists a positive constant $E_0$, depending on $N$, $m$, $\beta$, $R_0$ and $\delta_0$ only, such that
\begin{equation}
\mathcal{H}^{N-1}\left(\{\hat{x}_1\in\mathbb{S}^{N-1}:\ 
P_d\in S^+(\hat{x}_1)
\}\right)\geq E_0\min\{d^{(N-1)/2},1\}.
\end{equation}
\end{lem}

\proof{.} We sketch the proof of the lemma. First of all we notice that the distance of $P_d$ from $Q$ is bounded from below by $r_0+\cos(\theta_0)d$ and from above by $r_0+d$.
Let us call $T^+$ the portion of $\partial B_{\rho_0}(Q)$ that is formed by points $x\in \partial B_{\rho_0}(Q)$ such that the segment connecting $x$ to $P_d$ intersects $\partial B_{\rho_0}(Q)$ only at $x$. A simple computation shows that
$$\mathcal{H}^{N-1}(T^+)
\geq E_1\min\{d^{(N-1)/2},1\}$$
where $E_1$ is a positive constant depending on $N$, $m$, $\beta$, $R_0$ and $\delta_0$ only.
In fact, $T^+$ is the intersection of $\partial B_{\rho_0}(Q)$ with a symmetric cone with vertex in $Q$ and bisecting line passing through $P_d$ and whose amplitude is given by an angle $\alpha$,
$\alpha$ being of the order of $\sqrt{d}$

Let $\pi(\hat{x})$ be the supporting hyperplane at $P$. Let us call $D$ the open region which is enclosed
by  $\partial B_{\rho_0}(Q)$ and all tangent lines to $\partial B_{\rho_0}(Q)$ passing through $P_d$.
Notice that this is a portion of a symmetric cone with vertex in $P_d$ and bisecting line $l$ containing $Q$.
We wish to prove that
\begin{equation}\label{intermediatestep}
\mathcal{H}^{N-1}(\pi(\hat{x})\cap D)
\geq E_2\min\{d^{(N-1)/2},1\}
\end{equation}
where $E_2$ is a positive constant depending on $N$, $m$, $\beta$, $R_0$ and $\delta_0$ only.
In order to prove this property, let us begin with the following intermediate step.
We call $\Pi$ the plane of $\mathbb{R}^N$ containing $P$ and $Q$, and consequently $P_d$.
Let us take the two points $x_1$ and $x_2$ which are the intersections of $\pi(\hat{x})$ with the two lines in $\Pi$ passing through $P_d$ and tangent to $B_{\rho_0}(Q)\cap \Pi$. It is convenient to perform a rigid change of variables such that, in this new coordinate system,
$Q=0$ and $e_N=P_d-Q/\|P_d-Q\|$.
We show that
$$\|x_1-x_2\|\geq E_3\min\{\sqrt{d},1\}$$
where $E_3$ is a positive constant depending on $m$, $\beta$, $R_0$ and $\delta_0$ only.
In order to prove this, we begin with the case in which $P$ belongs to the segment connecting $Q$ to $P_d$. Then a simple geometric construction, using the properties of the angle $\alpha$ defined before, implies that $\|x_1-x_2\|=2\|x_1-P\|$ and $a_1=\|x_1-P\|$ is of order $\sqrt{d}$.
In the general case, we always have that $\|x_1-x_2\|\geq a_1$.
Let us notice that this concludes the proof of \eqref{intermediatestep} at least for $N=2$.

For $N>2$, the key step, which follows from elementary calculations, is to prove that the distance of $P$ from the line $l$ passing through $Q$ and $P_d$ is bounded by a constant times $d$. Then we take the point $x_3$ which is the intersection of the segment connecting $x_1$ to $x_2$ with $l$.
Then we construct the point $x_4$ which is one of the intersections of $\partial D\cap \Pi$ with
the hyperplane passing through $x_3$ with normal $P_d-Q$. Another computation leads to show that $\|x_3-x_4\|$ is of order $\sqrt{d}$. From this last property \eqref{intermediatestep} easily follows.

For any $\hat{x}_1\in\mathbb{S}^{N-1}$, let us call $l(\hat{x}_1)=\{x=s\hat{x}_1:\ s\geq 0\}$.
We have that if $l(\hat{x}_1)$ intersects $\pi(\hat{x})\cap D$ then $P_d\in S^+(\hat{x}_1)$.
Then the thesis immediately follows from \eqref{intermediatestep}.\cvd

\bigskip

We now finally restrict ourselves to $N=2,3$, $m\geq 2$ and $k_0\geq 1$ such that
$k_0R_0\geq z_0\geq 4$.
We let $k\geq k_0$ and, first of all, by \eqref{l2estimate}, we can estimate $\tilde{M}$ as follows
\begin{equation}\label{tildeMest}
\tilde{M}\leq C_5k
\end{equation}
where $C_5$ depends on $N$, $m$, $\beta$, $R_0$ and $\delta_0$ only.

Then we proceed in the following manner.
For any direction $\hat{x}\in\mathbb{S}^{N-1}$, the previous Lemma~\ref{hyperplanelemma} allows us to estimate
$$\int_{B_{\rho_1}\cap S^+(\hat{x})}|v|^2=
\int_{S^+(\hat{x})}|v|^2(x)\chi_{B_{\rho_1}}(x)\rmd x\leq C^2_4\eta^2_2(\varepsilon,k,\tilde{M})
$$
where $\chi$ denotes characteristic functions.
Therefore
$$
\int_{\mathbb{S}^{N-1}}
\left(\int_{S^+(\hat{x})}|v|^2(x)\chi_{B_{\rho_1}}(x)\rmd x\right)\rmd \mathcal{H}^{N-1}(\hat{x})
\leq \mathcal{H}^{N-1}(\mathbb{S}^{N-1})C_4^2\eta_2^2(\varepsilon,k,\tilde{M}).
$$
But, by Fubini's Theorem,
\begin{multline*}
\int_{\mathbb{S}^{N-1}}\!\!
\bigg(\int_{S^+(\hat{x})}|v|^2(x)\chi_{B_{\rho_1}}(x)\rmd x\bigg)\rmd \mathcal{H}^{N-1}(\hat{x})
=\\
\int_{\mathbb{S}^{N-1}}\!\!
\bigg(\int_{\mathbb{R}^N}|v|^2(x)\chi_{B_{\rho_1}}(x)\chi_{S^+(\hat{x})}(x)\rmd x\bigg)\rmd \mathcal{H}^{N-1}(\hat{x})=\\
\int_{\mathbb{R}^N} |v|^2(x)\chi_{B_{\rho_1}}(x)\bigg(\int_{\mathbb{S}^{N-1}}
\chi_{S^+(\hat{x})}(x)\rmd \mathcal{H}^{N-1}(\hat{x})
\bigg)\rmd x=\\
\int_{B_{\rho_1}} |v|^2(x)f(x)
\rmd x=\int_{B_{\rho_1}\backslash K} |v|^2(x)f(x)
\rmd x
\end{multline*}
where for any $x\in\mathbb{R}^N$
$$f(x)=\int_{\mathbb{S}^{N-1}}
\chi_{S^+(\hat{x})}(x)\rmd \mathcal{H}^{N-1}(\hat{x})$$
and we used the fact that $f(x)=0$ for any $x\in K$.

Then let us fix a constant $\mu$, $0<\mu\leq 1$.
By a repeated use of H\"older inequality we have
\begin{equation}
\int_{B_{\rho_1}\backslash K}|v|\leq
\left(\int_{B_{\rho_1}\backslash K}|v|^2f\right)^{a_1}
\left(\int_{B_{\rho_1}\backslash K}|v|^2\right)^{a_2}
\left(\int_{B_{\rho_1}\backslash K}f^{-\gamma}\right)^{1/2}
\end{equation}
where $0<\gamma\leq 1$, $a_1+a_2=1/2$ and their are given by the following formulas
\begin{equation}\label{coefficientsdefin}
\gamma=\frac{\mu}{2-\mu},\quad a_1=\frac{\mu}{2(2-\mu)},\quad a_2=\frac{1-\mu}{2-\mu}.
\end{equation}

The crucial remark is the following. For any $0<\mu<1$ we have that $\gamma=\mu/(2-\mu)<1$ and
there exists a constant $F(\mu)$, depending on $\mu$, $N$, $m$, $\beta$, $R_0$ and $\delta_0$ only, such that
\begin{equation}\label{crucialremarkf}
\int_{B_{\rho_1}\backslash K}f^{-\gamma}\leq F(\mu).
\end{equation}
We have that \eqref{crucialremarkf} follows from 
Lemma~\ref{Kproperties} and this construction. We integrate in spherical coordinates
\begin{multline}\label{crucialremarkf2}
\int_{B_{\rho_1}\backslash K}f^{-\gamma}=
\int_{\mathbb{S}^{N-1}}\left(
\int_{s_0(\hat{x})}^{\rho_1}f^{-\gamma}(s\hat{x})s^{N-1}\rmd s\right)\rmd \mathcal{H}^{N-1}(\hat{x})\leq\\
 \int_{\mathbb{S}^{N-1}}\left(
\int_{0}^{\rho_1-s_0(\hat{x})} (E_0\min\{s^{(N-1)/2},1\})^{-\gamma}(s_0(\hat{x})+s)^{N-1}\rmd s
\right)\rmd \mathcal{H}^{N-1}(\hat{x})\leq F(\mu).
\end{multline}
Let us note that, for $N=2$, $F(1)$ is also bounded  therefore we may allow $\mu=1$, hence $\gamma=1$, $a_1=1/2$, and
$a_2=0$.

Then, setting $C_6=C_4^2\mathcal{H}^{N-1}(\mathbb{S}^{N-1})$ and recalling Lemma~\ref{hyperplanelemma},
 we may conclude that
\begin{equation}\label{N=23est}
\int_{B_{\rho_1}\backslash K}|v|\leq\\
\left\{\begin{array}{ll}
F(1)^{1/2}C_6^{1/2}
\eta_2(\varepsilon,k,\tilde{M})
& N=2\\
&\\
F(\mu)^{1/2}C_6^{a_1}
\tilde{M}^{2a_2}
(\eta_2(\varepsilon,k,\tilde{M}))^{2a_1} & N=3,\ 0<\mu<1.
\end{array}\right.
\end{equation}

Thus we have proved the following stability result.
 
\begin{teo}\label{uptheboundaryteo}
Under the previous notation and assumptions, let us assume that
$N=2,3$, $m\geq 2$ and $k_0\geq 1$ such that
$k_0R_0\geq z_0\geq 4$.
We let $k\geq k_0$, $\omega\in\mathbb{S}^{N-1}$ and $\Sigma_1$ and
$\Sigma_2\in X(m,\beta,R_0,\delta_0)$.

Then there exist positive constants $C_1,\ldots,C_6$, depending on $N$, $m$, $\beta$, $R_0$ and $\delta_0$ only, such that the following holds.

Let us assume that
$$\|\mathcal{A}(\Sigma_1)(\cdot,\omega,k)-\mathcal{A}(\Sigma_2)(\cdot,\omega,k)\|_{L^2(\mathbb{S}^{N-1})}\leq\varepsilon\leq C_5k
$$
and that
$$\eta_1(\varepsilon,k,C_5k,\tilde{B}_1)\leq (C_3/C_2)C_5k.
$$

Then we have the stability estimate
\begin{equation}\label{stabestimate}
\|u_1-u_2\|_{L^1(B_{\rho_1}\backslash K)}\leq \eta\left(\|\mathcal{A}(\Sigma_1)(\cdot,\omega,k)-\mathcal{A}(\Sigma_2)(\cdot,\omega,k)\|_{L^2(\mathbb{S}^{N-1})}\right)
\end{equation}
where $\eta$ satisfies the following.
Let $\eta_2=\eta_2(\varepsilon,k,C_5k)$ be given as in Lemma~\textnormal{\ref{hyperplanelemma}}.

For $N=2$, we have
\begin{equation}\label{N=2finalest}
\eta(\varepsilon)\leq
F(1)^{1/2}C_6^{1/2}\eta_2(\varepsilon,k,C_5k).
\end{equation}

For $N=3$, fixed $\mu$, $0<\mu<1$, we have
\begin{equation}\label{N=3finalest}
\eta(\varepsilon)\leq
F(\mu)^{1/2}C_6^{a_1}
C_5^{2a_2}
k^{2a_2}
(\eta_2(\varepsilon,k,C_5k))^{2a_1},
\end{equation}
where $a_1$ and $a_2$ are given in \eqref{coefficientsdefin}.
\end{teo}

We conclude this section by making
the stability estimate of Theorem~\ref{uptheboundaryteo} more explicit at least in the extremely high frequencies regime. Namely, let us assume that
$N=2,3$, $m\geq 2$ and $k_0\geq 1$ such that
$k_0R_0\geq z_0\geq 4$.
We let $k\geq k_0$, $\omega\in\mathbb{S}^{N-1}$ and $\Sigma_1$ and
$\Sigma_2\in X(m,\beta,R_0,\delta_0)$.

We begin by noticing that \eqref{Massum} holds for any $k\geq k_0$, with $\tau=1$ and a constant 
depending on $N$, $m$, $\beta$, $R_0$ and $\delta_0$ only, namely
$$M\leq\tilde{M}\leq C_5 k,$$
$C_5$ as in the previous theorem. Hence,
if $\varepsilon\leq 1/\rme$ and
$$k\geq \frac{\tilde{C}_0}{\log(4/3)\tilde{B}_0R_0}\log(1/\varepsilon),$$
we can apply Corollary~\ref{Lipschitz} and obtain that
$$\eta_1(\varepsilon,k,C_5k)
\leq\left( (\rho_2-\rho_1)
(\tilde{A}\tilde{B}_1(R+1)(1+C(1)C^2_5 (\tilde{C}_0/\tilde{B}_0R_0)^{2})/R)\right)^{1/2}
\varepsilon = C_8\varepsilon,
$$
where $C(1)$ is an absolute constant and $C_8$ clearly depends on $N$, $m$, $\beta$, $R_0$ and $\delta_0$ only. 
Hence, if $\varepsilon\leq C_3C_5k/(C_2C_8)$, we have that
\begin{multline}\label{finalexample}
\|u_1-u_2\|_{L^1(B_{\rho_1}\backslash K)}\leq\\
F(\mu)^{1/2}
C_6^{a_1}
C_5^{2a_2}
k^{2a_2}
\left((C_2C_8)^2k^2\varepsilon^2
+\frac{
(C_3C_5)^2k^4}{\left(-\log\left(C_2C_8\varepsilon/(C_3C_5k)\right)+k\right)^{1/8}}
\right)^{a_1}
\end{multline}
where $0<\mu< 1$ for $N=3$ and $\mu=1$ for $N=2$.

We notice that we have an estimate with an explicit dependence on $k$, in the extremely high frequencies regime. However we have to point out that we do not have any increasing stability phenomenon as $k$ grows.

\section{Instability for the inverse scattering problem}\label{instabilitysec}

In this section we fix $N=2,3$. We also fix an integer $m\geq 2$ and positive constants $\beta$ and $R_0$. We fix $\delta_0>0$ as defined in Proposition~\ref{discretesetprop} and we set $R=R_0+\delta_0$.
We fix a positive constant $k_0$ and we denote
\begin{equation}\label{I_Ntildedefin}
\tilde{I}_N=\left\{\begin{array}{ll}
[k_0,+\infty)&\text{if }N=2,\\
(0,+\infty)&\text{if }N = 3.
\end{array}\right.
\end{equation}

We fix $\Sigma\in X(m,\beta,R_0,\delta_0)$ and we consider its far-field pattern $\mathcal{A}(\Sigma)$ and its decomposition in spherical harmonics. Then, by \eqref{decompscatcoeff} and \eqref{link} we have, for any $\omega\in\mathbb{S}^{N-1}$, any $k>0$, any index $i$,
and any $r>R$
$$|\tilde{b}_i(\omega,k)|\leq \sqrt{\frac{2}{\pi}}k^{-(N-1)/2}\frac{(kr)^{(N-2)/2}}{|H^{(1)}_{\gamma(v_i)+(N-2)/2}(kr)|}\left|\int_{\mathbb{S}^{N-1}}u^s(r\hat{x};\omega,k,\Sigma)
v_i(\hat{x})\mathrm{d}\hat{x}\right|.
$$
Hence, by collecting the estimates of the previous section, we have
for any $\omega\in\mathbb{S}^{N-1}$, any $k\in\tilde{I}_N$, and any index $i$
$$|\tilde{b}_i(\omega,k)|\leq C_1(kr)^{-1/2}\frac{1}{|H^{(1)}_{\gamma(v_i)+(N-2)/2}(kr)|}
\max\{1,k^{(N-1)/2}\}
,\quad\text{for any }r\geq R+1,
$$
where $C_1$ depends on $N$, $m$, $\beta$, $R_0$, $R$ and, only if $N=2$, $k_0$.
We call $\tilde{R}=R+1$.
If $k\leq 1/\tilde{R}$ we choose $r=1/k$ and we obtain
$$|\tilde{b}_i(\omega,k)|\leq \frac{C_1}{|H^{(1)}_{\gamma(v_i)+(N-2)/2}(1)|},\quad \text{for }k\in\tilde{I}_N,\ k\leq 1/\tilde{R}.
$$
If $k\geq 1/\tilde{R}$ we choose $r=\tilde{R}$ and we obtain
$$|\tilde{b}_i(\omega,k)|\leq C_1\frac{(k\tilde{R})^{(N-2)/2}}{|H^{(1)}_{\gamma(v_i)+(N-2)/2}(k\tilde{R})|},\quad \text{for }k\in\tilde{I}_N,\ k\geq 1/\tilde{R}.
$$
Finally, if we set $z(k)=\max\{1,k\tilde{R}\}$, we have for any $\omega\in\mathbb{S}^{N-1}$ and any index $i$
\begin{equation}
|\tilde{b}_i(\omega,k)|\leq C_1\frac{(z(k))^{(N-2)/2}}{|H^{(1)}_{\gamma(v_i)+(N-2)/2}(z(k))|}
,\quad \text{for }k\in \tilde{I}_N.
\end{equation}
On the other hand, we also have for any $\omega\in\mathbb{S}^{N-1}$ and any index $i$
\begin{equation}\label{uniformL2bound1}
|\tilde{b}_i(\omega,k)|\leq C_2,\quad \text{for }k\in \tilde{I}_N
\end{equation}
where $C_2$ depends on $N$, $m$, $\beta$, $R_0$, $R$ and, only if $N=2$, $k_0$.

Then we use Corollary~\ref{unifbounded} and Theorem~\ref{regimes} to obtain the following result.
For any $k\in\tilde{I}_N$, and any arbitrary $\omega\in\mathbb{S}^{N-1}$, we have two different cases.
For any index $i$ such that $\gamma(v_i)\geq \rme z(k)$ we have
\begin{equation}
|\tilde{b}_i(\omega,k)|\leq C_3
(z(k))^{(N-1)/2}
\left(\frac{a\nu(i)}{\rme z(k)}\right)^{-(\nu(i)-1/2)}
\end{equation}
where $\nu(i)=\gamma(v_i)+(N-2)/2$, $a=1+\sqrt{\rme^2-1}/\rme$ and $C_3$ depends on $N$, $m$, $\beta$, $R_0$, $R$ and, only if $N=2$, $k_0$. Let us also notice that $\nu(i)\geq \rme$, therefore $\nu(i)-1/2\geq 2$ and, obviously, $a\nu(i)/(\rme z(k))\geq a>1$.

We note that there exists a constant $\tilde{c}\geq 1$, depending on $N$ only, such that for any $z\geq 1$ and any $t$ such that $t\geq
\tilde{c}\rme z$ we have 
$$z^{(N-1)/2}
\left(a(t+(N-2)/2)/(\rme z)\right)^{-(t+(N-2)/2-1/2)}\leq
z^{(N-1)/2}
\left(at/(\rme z)\right)^{-(t+(N-2)/2-1/2)}
\leq 1.$$
We notice that, since $\rme\log(a)>1$, for $N=2,3$ we may actually choose $\tilde{c}=1$.
On the other hand, we recall that
for any index $i$ we have \eqref{uniformL2bound1}.

Let $\tilde{C}= (\mathcal{H}^{N-1}(\mathbb{S}^{N-1})^{1/2}\max\{C_2,C_3\}$.
Obviously, $\tilde{C}$ depends on $N$, $m$, $\beta$, $R_0$, $R$ and, only if $N=2$, $k_0$. 
Without loss of generality we may assume that $\tilde{C}\geq 2$.

By the reciprocity relation \eqref{recrel} we conclude that for any $k\in\tilde{I}_N$, and for any indexes $i,l$, we have
\begin{equation}
|b_{i,l}(k)|\leq \tilde{C}
\end{equation}
and
for any $k\in\tilde{I}_N$, and for any indexes $i,l$ such that $\max\{\gamma(v_i),\gamma(v_l)\}\geq \tilde{c}\rme z(k)$, we have
\begin{equation}
|b_{i,l}(k)|\leq \tilde{C}
(z(k))^{(N-1)/2}
\left(\frac{a\max\{\gamma(v_i),\gamma(v_l)\}}{\rme z(k)}\right)^{-(\max\{\gamma(v_i),\gamma(v_l)\}+(N-3)/2)}\leq \tilde{C}.
\end{equation}

For any $\Sigma\in X(m,\beta,R_0,\delta_0)$ and any $k\in\tilde{I}_N$, we have $\mathcal{A}(\Sigma)(\cdot,\cdot,k)\in Y_s(\mathbb{S}^{N-1}\times \mathbb{S}^{N-1})$ for any $s\geq 0$.
We also recall that, again for any $s\geq 0$,
$$\|\mathcal{A}(\Sigma)(\cdot,\cdot,k)\|_{H^s}=\|\mathcal{A}(\Sigma)(\cdot,\cdot,k)\|_{H^s(\mathbb{S}^{N-1}\times \mathbb{S}^{N-1})}\leq C_4\|\mathcal{A}(\Sigma)(\cdot,\cdot,k) \|_s$$
where $C_4=4$.

For any fixed $k\in\tilde{I}_N$,
we denote the set
$$\tilde{Y}(k)=\{\mathcal{A}(\Sigma)(\cdot,\cdot,k):\ \Sigma\in
X(m,\beta,R_0,\delta_0)\}.$$
We notice that for any fixed $s\geq 0$, $\tilde{Y}(k)
\subset H^s(\mathbb{S}^{N-1}\times \mathbb{S}^{N-1})$ and it may be considered as a metric space endowed with the distance induced by the $H^s$ norm.
We recall that, for a given positive $\varepsilon$, 
a subset $\tilde{Y}'\subset \tilde{Y}(k)$ is an $\varepsilon$-\emph{net} for 
$\tilde{Y}(k)$, with respect to the $H^s$ norm, if for every $y\in \tilde{Y}(k)$ there exists $y'\in \tilde{Y}'$
whose $H^s$ distance from $y$ is less than or equal to $\varepsilon$.

Before stating the main instability theorem, let us introduce the following notation.
We begin by noticing that for any $t$ such that $t\geq \max\{\tilde{c}\rme z(k),2s+N\}$
 we have that
$$f(t)=(1+t)^{2s+N-1/2}\left(\frac{at}{\rme z(k)}\right)^{-(t+(N-3)/2)}$$
is a decreasing function of $t$. We call
\begin{equation}\label{epsk}
\tilde{\varepsilon}(k)=2C_4\tilde{C}
(z(k))^{(N-1)/2}f(\max\{\tilde{c}\rme^2 ,4s+(3N/2)+1\}z(k)).
\end{equation}
Let us notice that $\tilde{\varepsilon}(k)$ is a decreasing function of $k$ which goes to $0$ as $k\to\infty$. 
Finally, we denote
\begin{equation}\label{Zkdef}
\begin{array}{c}
\tilde{B}(s)=\max\{\tilde{c}\rme^2 ,4s+(3N/2)+1\}\quad
Z(k)=\tilde{B}(s)\max\{1,(R+1)k\}\\ C_5=(2C_4)(2\tilde{C}+1).
\end{array}
\end{equation}
We observe that $\tilde{B}(s)$ depends on $s$ and $N$ only, whereas $C_5$ depends on 
$N$, $m$, $\beta$, $R_0$, $R$ and, only if $N=2$, $k_0$.

\begin{teo}\label{instmainteo}
Fixed $s\geq 0$ and $k\in \tilde{I}_N$,
for every $\varepsilon$, $0<\varepsilon<1/\rme$, there exists $\delta=\delta(\varepsilon,k)$, $0<\delta\leq \delta_0$, and two obstacles $\Sigma_1$ and $\Sigma_2$ belonging to $X(m,\beta,R_0,\delta)$ such that
\begin{equation}\label{finalinst}
d_H(\Sigma_1,\Sigma_2)\geq \delta\quad\text{and}\quad
\|\mathcal{A}(\Sigma_1)(\cdot,\cdot,k)-\mathcal{A}(\Sigma_2)(\cdot,\cdot,k)\|_{H^s(\mathbb{S}^{N-1}\times \mathbb{S}^{N-1})}\leq 2\varepsilon.
\end{equation}

If $\varepsilon\geq \tilde{\varepsilon}(k)$, then
\begin{equation}\label{deltadefepslarge}
\delta(\varepsilon,k)=\delta_0\frac{2^{-m(N+3)/(N-1)}}{(1+Z(k))^{2m}}
\left[\log\left(C_5(1+Z(k))^{(2s+N-1/2)}
/\varepsilon\right)\right]^{-m/(N-1)},
\end{equation}
hence
\begin{multline}\label{finalinstestepslarge}
d_H(\Sigma_1,\Sigma_2)\geq \\
\delta_0\frac{2^{-m(N+3)/(N-1)}}{(1+Z(k))^{2m}}\left[\log\left(\frac{2C_5(1+Z(k))^{(2s+N-1/2)}
}{\|\mathcal{A}(\Sigma_1)(\cdot,\cdot,k)-\mathcal{A}(\Sigma_2)(\cdot,\cdot,k)\|_{H^s}}\right)\right]^{-m/(N-1)}.
\end{multline}

If $0<\varepsilon< \tilde{\varepsilon}(k)$, then
\begin{equation}\label{deltadefepssmall}
\delta(\varepsilon,k)=\delta_0\frac{2^{-m(N+3)/(N-1)}}{(1+\tilde{t})^{2m}}\left[\log\left(C_5(1+\tilde{t})^{(2s+N-1/2)}
/\varepsilon\right)\right]^{-m/(N-1)}
\end{equation}
hence
\begin{multline}\label{finalinstestepssmall}
d_H(\Sigma_1,\Sigma_2)\geq \\
\delta_0\frac{2^{-m(N+3)/(N-1)}}{(1+\tilde{t})^{2m}}\left[\log\left(\frac{2C_5(1+\tilde{t})^{(2s+N-1/2)}
}{\|\mathcal{A}(\Sigma_1)(\cdot,\cdot,k)-\mathcal{A}(\Sigma_2)(\cdot,\cdot,k)\|_{H^s}}\right)\right]^{-m/(N-1)},
\end{multline}
where $\tilde{t}> Z(k)$ satisfies
\begin{equation}\label{deltatdef}
2C_4\tilde{C}
(z(k))^{(N-1)/2}f(\tilde{t})=\varepsilon.
\end{equation}
\end{teo}

\proof{.}
Let us then fix $s\geq 0$, $k\in \tilde{I}_N$, and
$\varepsilon$, $0<\varepsilon<1/\rme$. The crucial step is constructing 
an $\varepsilon$-net
for $\tilde{Y}(k)$, with respect to the $H^s$ norm, and estimating its number of elements.
We distinguish two regimes. First we treat the case when $\varepsilon\geq \tilde{\varepsilon}(k)$, then we shall deal with the
case $0<\varepsilon<\tilde{\varepsilon}(k)$.

If $\varepsilon\geq \tilde{\varepsilon}(k)$, for any integer $n$ such that $n\geq Z(k)$ we have
$$2C_4\tilde{C}
(z(k))^{(N-1)/2}f(n)\leq\tilde{\varepsilon}(k)\leq  \varepsilon.
$$

Let $\tilde{n}$ be the integer part of $Z(k)$.
Let 
$\varepsilon'=(1+\tilde{n})^{-(2s+N-1/2)}\varepsilon/(2C_4)$
and $\Psi_{\varepsilon}=[-\tilde{C},\tilde{C}]\cap\varepsilon'\mathbb{Z}$. We remark that $\Psi_{\varepsilon}$ is a
finite subset of $\mathbb{R}$ and we have that $\#\Psi_{\varepsilon}\leq (2\tilde{C}+1)/\varepsilon'$.

Let us define the following subset of $L^2(\mathbb{S}^{N-1}\times \mathbb{S}^{N-1})$
$$\hat{Y}(\varepsilon)=\{g\in L^2(\mathbb{S}^{N-1}\times \mathbb{S}^{N-1}):\ a_{i,l}\in\Psi_{\varepsilon}\text{ if }
\max\{\gamma(v_i),\gamma(v_l)\}\leq\tilde{n}\text{ and }a_{i,l}=0\text{ otherwise}\}.$$
We may count the number of elements of $\hat{Y}(\varepsilon)$ as follows. If we set
$$s=\#\{(i,l):\ \max\{\gamma(v_i),\gamma(v_l)\}\leq\tilde{n}\}$$
we obtain that
$$s\leq 4(1+\tilde{n})^{2N-2}.$$
Then we have that
$\#\hat{Y}(\varepsilon)=(\#\Psi_{\varepsilon})^{s}$ and hence
\begin{multline*}
\#\hat{Y}(\varepsilon)\leq ((2\tilde{C}+1)/\varepsilon')^s\leq
\left((2C_4)(2\tilde{C}+1)(1+\tilde{n})^{(2s+N-1/2)}
/\varepsilon\right)^s\leq \\
\left(C_5(1+\tilde{n})^{(2s+N-1/2)}
/\varepsilon\right)^{4(1+\tilde{n})^{2N-2}}.
\end{multline*}

It is now easy to construct an $\varepsilon$-net for $\tilde{Y}(k)$ with respect to the $H^s$ norm with at most 
$$\exp\left({4(1+Z(k))^{2N-2}}\log\left(C_5(1+Z(k))^{(2s+N-1/2)}
/\varepsilon\right)\right)
$$ elements, see for instance the proof of Lemma~2.3 in \cite{DC-R1}.

We conclude the proof of the theorem in the first regime.
Let us assume that there exists $\delta$, $0<\delta\leq \delta_0$, such that
$$\exp(2^{-N}(\delta_0/\delta)^{(N-1)/m})> 
\exp\left(4(1+Z(k))^{2N-2}\log\left(C_5(1+Z(k))^{(2s+N-1/2)}
/\varepsilon\right)\right).
$$
Then, there exist two elements $\Sigma_1$ and $\Sigma_2$ of $X(m,\beta,R_0,\delta)$
satisfying \eqref{finalinst}.
This is true if
$$2^{-N}(\delta_0/\delta)^{(N-1)/m}>
4(1+Z(k))^{2N-2}\log\left(C_5(1+Z(k))^{(2s+N-1/2)}
/\varepsilon\right)$$
that is, for instance, when $\delta$ is given by \eqref{deltadefepslarge}. Therefore, also \eqref{finalinstestepslarge} immediately follows and the theorem is fully proved in the first regime.

Let us now consider the second regime, that is when 
$0<\varepsilon<\tilde{\varepsilon}(k)$. Let $\tilde{t}>Z(k)$ be such that
$2C_4\tilde{C}
(z(k))^{(N-1)/2}f(\tilde{t})=\varepsilon$.
Then we repeat exactly the same procedure just by replacing $\tilde{n}$ with the
integer part of $\tilde{t}$ and the proof is concluded.\cvd

\bigskip

We conclude this section with a few comments on the results contained in the previous instability theorem. First of all we make estimates \eqref{deltadefepssmall} and \eqref{finalinstestepssmall}
more readable by estimating in a suitable way $\tilde{t}$.
In the theorem, we are assuming the wavenumber $k$ fixed and we are establishing how the instability changes with respect to the error $\varepsilon$. However, in order to understand the high frequency asymptotics, we then consider $\varepsilon$ to be fixed and discuss the changes in the instability as $k$ increases.

We being with the first remark.
We notice
that, since $\varepsilon<\tilde{\varepsilon}(k)$, we have
$\tilde{t}\geq Z(k)$ and consequently
$$f(\tilde{t})\leq f_1(\tilde{t})=\left(\frac{3\rme z(k)}{2a}\right)^{2s+N-1/2}\left(\frac{a\tilde{t}}{\rme z(k)}\right)^{-\tilde{t}/2}
.$$
It is enough to
find $\hat{t}$ such that
$$2C_4\tilde{C}
(z(k))^{(N-1)/2}f_1(\hat{t})= \varepsilon$$
to deduce that
$Z(k)\leq \tilde{t}\leq \hat{t}$. Hence our result holds true if we replace $\tilde{t}$ with $\hat{t}$ in \eqref{deltadefepssmall} and \eqref{finalinstestepssmall}. Finally, a straightforward computation shows that
$$\hat{t}\leq 2\log\left(\frac{\tilde{b}(k,s)}{\varepsilon}\right)$$
where
\begin{equation}\label{tildebdefin}
\tilde{b}(k,s)=2C_4\tilde{C}(3\rme/(2a))^{2s+N-1/2}
(z(k))^{2s+(3N/2)-1}.
\end{equation}

Let us now notice that
\begin{equation}\label{intermediateeq}
C_5(1+\tilde{t})^{2s+N-1/2}/\varepsilon=
\frac{2\tilde{C}+1}{\tilde{C}(z(k))^{(N-1)/2}}
\left(\frac{a\tilde{t}}{\rme z(k)}\right)^{\tilde{t}+(N-3)/2}.
\end{equation}
We deduce that
$$C_5(1+\tilde{t})^{2s+N-1/2}/\varepsilon\leq 
\frac{2\tilde{C}+1}{\tilde{C}(z(k))^{(N-1)/2}}
\left(\frac{a\hat{t}}{\rme z(k)}\right)^{(8/3)(\hat{t}/2)}=
\frac{2\tilde{C}+1}{\tilde{C}(z(k))^{(N-1)/2}}\left(\frac{\tilde{b}(k,s)}{\varepsilon}\right)^{8/3}.
$$

Since $\tilde{C}\geq 2$ and $z(k)\geq 1$ we conclude that
we can replace $\delta(\varepsilon,k)$ in \eqref{deltadefepssmall} by
$$\delta=
\frac{\delta_0}{2^{m(N+3)/(N-1)}}\left(1+(8/3)\log\left(\frac{\tilde{b}(k,s)}{\varepsilon}\right)\right)^{-2m-m/(N-1)}
$$
and \eqref{finalinstestepssmall} may be replaced by
\begin{multline*}
d_H(\Sigma_1,\Sigma_2)\geq \\
\frac{\delta_0}{2^{m(N+3)/(N-1)}}\left(1+(8/3)\log\left(\frac{2\tilde{b}(k,s)}{\|\mathcal{A}(\Sigma_1)(\cdot,\cdot,k)-\mathcal{A}(\Sigma_2)(\cdot,\cdot,k)\|_{H^s}}
\right)\right)^{-2m-m/(N-1)}
\end{multline*}
where $\tilde{b}(k,s)$ is as in \eqref{tildebdefin}.

Let us now consider a fixed $\varepsilon$, $0<\varepsilon<1/\rme$. Let us consider that for some $k\in \tilde{I}_N$ we have $\tilde{\varepsilon}(k)>\varepsilon$, otherwise we have that
$\varepsilon\geq \tilde{\varepsilon}(k)$ for any $k\in\tilde{I}_N$ and it is easy to understand the instability behavior since we are always in the first regime and \eqref{deltadefepslarge} 
and \eqref{finalinstestepslarge} apply.

Then let $k(\varepsilon)\in \tilde{I}_N$ be the first $k\in \tilde{I}_N$ such that $\tilde{\varepsilon}(k(\varepsilon))=\varepsilon$. An easy computation shows that $k(\varepsilon)$ grows
essentially like a constant times $\log(1/\varepsilon)$.
As long as $k<  k(\varepsilon)$ we have that the improvement in the instability as $k$ increases is not very big since we need to use \eqref{deltadefepssmall} and \eqref{finalinstestepssmall}. Let us notice that as
$k<  k(\varepsilon)$ we have $Z(k)<\tilde{t}\leq Z(k(\varepsilon))$, therefore as
$k<k(\varepsilon)$ increases and converges to 
$k(\varepsilon)$ we have that $\tilde{t}$ increases and converges to $Z(k(\varepsilon))$.
From $k(\varepsilon)$ onwards, that is in the very high frequencies regime, 
\eqref{deltadefepslarge} 
and \eqref{finalinstestepslarge} apply and the improvement in the instability is more evident.
We state these observations in this final corollary.

\begin{cor}\label{instcorollary}
Let us fix $\varepsilon$, $0<\varepsilon<1/\rme$.
Let us assume that for some $k\in \tilde{I}_N$ we have $\tilde{\varepsilon}(k)>\varepsilon$ and
let $k(\varepsilon)\in \tilde{I}_N$ be the first $k\in \tilde{I}_N$ such that $\tilde{\varepsilon}(k(\varepsilon))=\varepsilon$.

If $k<  k(\varepsilon)$, then
$$
\delta=
\frac{\delta_0}{2^{m(N+3)/(N-1)}}\left(1+(8/3)\log\left(\frac{\tilde{b}(k,s)}{\varepsilon}\right)\right)^{-2m-m/(N-1)}
$$
hence
\begin{multline*}
d_H(\Sigma_1,\Sigma_2)\geq \\
\frac{\delta_0}{2^{m(N+3)/(N-1)}}\left(1+(8/3)\log\left(\frac{2\tilde{b}(k,s)}{\|\mathcal{A}(\Sigma_1)(\cdot,\cdot,k)-\mathcal{A}(\Sigma_2)(\cdot,\cdot,k)\|_{H^s}}
\right)\right)^{-2m-m/(N-1)}
\end{multline*}
where $\tilde{b}(k,s)$ is as in \eqref{tildebdefin}.

If $k\geq k(\varepsilon)$, then
$$
\delta=\delta_0\frac{2^{-m(N+3)/(N-1)}}{(1+Z(k))^{2m}}
\left[\log\left(C_5(1+Z(k))^{(2s+N-1/2)}
/\varepsilon\right)\right]^{-m/(N-1)},
$$
hence
\begin{multline*}
d_H(\Sigma_1,\Sigma_2)\geq \\
\delta_0\frac{2^{-m(N+3)/(N-1)}}{(1+Z(k))^{2m}}\left[\log\left(\frac{2C_5(1+Z(k))^{(2s+N-1/2)}
}{\|\mathcal{A}(\Sigma_1)(\cdot,\cdot,k)-\mathcal{A}(\Sigma_2)(\cdot,\cdot,k)\|_{H^s}}\right)\right]^{-m/(N-1)}.
\end{multline*}
\end{cor}

\end{document}